\documentclass[12pt]{article}

\usepackage[T1]{fontenc}
\usepackage{amsmath,amsthm,mathtools}

\IfFileExists{newtxtext.sty}{%
  \usepackage{newtxtext}
  \usepackage{newtxmath}%
}{%
  \usepackage{mathptmx}\usepackage{amssymb}%
}
\usepackage[dvipsnames]{xcolor}
\usepackage{tikz} 

\usepackage{iftex}
\ifPDFTeX\usepackage[utf8]{inputenc}\fi
\usepackage[a4paper,margin=1in]{geometry}
\usepackage{aliascnt}
\usepackage{enumitem}
\usepackage{needspace}
\usepackage{array}
\usepackage{booktabs}
\usepackage[protrusion=true,expansion=false]{microtype}
\usepackage[hyperfootnotes=false,colorlinks=true,linkcolor=blue,citecolor=blue,
  urlcolor=blue]{hyperref}

\makeatletter
\newcommand\blfootnote[1]{%
  \begingroup
  \renewcommand\@makefntext[1]{\noindent\strut ##1}%
  \renewcommand\thefootnote{}\footnote{#1}%
  \addtocounter{footnote}{-1}%
  \endgroup
}
\long\def\@makecaption#1#2{%
  \vskip\abovecaptionskip
  \sbox\@tempboxa{#1. #2}%
  \ifdim \wd\@tempboxa >\hsize
    #1. #2\par
  \else
    \global \@minipagefalse
    \hb@xt@\hsize{\hfil\box\@tempboxa\hfil}%
  \fi
  \vskip\belowcaptionskip}
\makeatother

\setlist{itemsep=0.25em,topsep=0.4em}
\newtheorem{theorem}{Theorem}[section]
\newaliascnt{proposition}{theorem}
\newtheorem{proposition}[proposition]{Proposition}
\aliascntresetthe{proposition}
\newaliascnt{lemma}{theorem}
\newtheorem{lemma}[lemma]{Lemma}
\aliascntresetthe{lemma}
\newaliascnt{corollary}{theorem}
\newtheorem{corollary}[corollary]{Corollary}
\aliascntresetthe{corollary}
\newaliascnt{question}{theorem}
\newtheorem{question}[question]{Question}
\aliascntresetthe{question}
\newaliascnt{remark}{theorem}
\newtheorem{remark}[remark]{Remark}
\aliascntresetthe{remark}
\newaliascnt{problem}{theorem}

\aliascntresetthe{problem}

\usepackage[nameinlink,capitalise,noabbrev]{cleveref}

\newcommand{\R}{\mathbb R}
\newcommand{\Z}{\mathbb Z}
\newcommand{\Q}{\mathbb Q}
\newcommand{\C}{\mathbb C}
\newcommand{\F}{\mathbb F}
\newcommand{\cU}{\mathcal U}
\newcommand{\cV}{\mathcal V}
\newcommand{\cW}{\mathcal W}
\newcommand{\cX}{\mathcal X}
\newcommand{\cE}{\mathcal E}
\newcommand{\hn}{\bar h}
\newcommand{\lcyl}{\ell_{\mathrm{cyl}}}
\newcommand{\ip}[2]{\langle #1,#2\rangle}
\newcommand{\Sc}{\mathcal R}
\newcommand{\cB}{\mathcal B}
\newcommand{\sign}{\operatorname{sign}}
\newcommand{\defic}{\operatorname{def}}
\newcommand{\step}[1]{\par\medskip\noindent\emph{#1}}
\numberwithin{equation}{section}
\hypersetup{pdftitle={Fundamental groups of closed four-manifolds with positive biorthogonal curvature},pdfauthor={Zhiqi Chen and Hui Zhang}}
\title{\textbf{Fundamental groups of closed four-manifolds with positive biorthogonal curvature}}
\author{Zhiqi Chen and Hui Zhang}
\date{}
\begin{document}
\raggedbottom
\clubpenalty=10000
\widowpenalty=10000
\displaywidowpenalty=10000
\maketitle
\begin{abstract}
In this paper, we prove that a group is the fundamental group of a closed, connected,
oriented smooth four-manifold admitting a metric of positive biorthogonal
curvature if and only if it is finitely presented. This answers a question
of Bettiol. The manifold can be chosen to be stably parallelizable, hence
spin with vanishing signature, and with Euler characteristic $2(1-p+q)$,
where $p$ and $q$ are the numbers of generators and of nontrivial relators
in a given finite presentation. Moreover, we show that the class of closed
oriented four-manifolds admitting positive biorthogonal curvature is closed
under surgery along embedded circles. As applications, we obtain metrics of
positive biorthogonal curvature on the boundary of every compact oriented
five-manifold that has a handle decomposition with handles of index at most
two, which includes every oriented $S^2$-bundle over a closed oriented
surface. We also obtain such metrics on $Y\times S^1$ for every closed
orientable three-manifold $Y$ of positive scalar curvature and, after
connected sum with copies of $S^2\times S^2$, on a large class of
four-manifolds with arbitrary fundamental group. Finally, we show that for
the fundamental group of a finite aspherical two-dimensional complex $X$
the minimal Euler characteristic of such manifolds equals $2\chi(X)$.
\end{abstract}
\blfootnote{%
}

\section{Introduction}\label{sec:introduction}

Which groups occur as fundamental groups of closed manifolds satisfying a
given curvature condition is a classical question, and in dimension four
the answer is known at both ends of the scale of positivity conditions.
For positive scalar curvature, Gromov and Lawson~\cite{GL80} and Schoen and
Yau~\cite{SY79} proved that the condition is preserved under surgery of
codimension at least three, and Carr~\cite[Corollary~2]{Carr88} proved
that every finitely presented group is the fundamental group of a closed
four-manifold of positive scalar curvature. Carr's manifolds are
boundaries of regular neighborhoods of finite two-dimensional complexes in
$\R^5$. They can also be obtained from a connected sum of copies of
$S^1\times S^3$ by surgery along circles that represent the relators of a
presentation. Since the fundamental group of a closed manifold is always
finitely presented, positive scalar curvature imposes no restriction at
all. Much stronger conditions, by contrast, leave very little room. For
positive isotropic curvature, introduced by Micallef and
Moore~\cite{MM88}, Fraser~\cite{Fraser03} proved that in dimensions at
least five the fundamental group of a closed manifold contains no subgroup
isomorphic to $\Z\oplus\Z$, and in dimension four the classification of
Chen, Tang and Zhu~\cite[Corollary~1]{CTZ} implies that the fundamental
group is \emph{virtually free}, that is, it contains a free subgroup of
finite index.

In this paper we determine the fundamental groups for a condition that
lies strictly between these two. Let $(M,g)$ be a Riemannian four-manifold.
The \emph{biorthogonal curvature} of a two-plane $\sigma\subset T_pM$ is
the average
\[
 K_g^\perp(\sigma)=\frac12\bigl(K_g(\sigma)+K_g(\sigma^{\perp})\bigr),
 \qquad\text{and we put}\qquad
 k_g(p)=\min_{\sigma\subset T_pM}K_g^\perp(\sigma),
\]
where $K_g$ denotes the sectional curvature and $\sigma^{\perp}$ the
orthogonal complement of $\sigma$ in $T_pM$. We say that $g$ has
\emph{positive biorthogonal curvature} if $k_g>0$ everywhere. In an
orthonormal frame, the scalar curvature $s_g$ is four times the sum of the
biorthogonal curvatures of the three pairs of complementary coordinate
planes, and hence
\[
 K_g>0\quad\Longrightarrow\quad k_g>0\quad\Longrightarrow\quad s_g>0.
\]
Neither implication can be reversed. Biorthogonal curvature was studied
by Seaman~\cite{Seaman91} in connection with orthogonal pinching, and by
Costa~\cite{Costa12} and Costa--Ribeiro~\cite{CR14} in connection with
modified Yamabe invariants. On the existence side,
Bettiol~\cite{Bettiol14} constructed metrics with $k_g>0$ on
$S^2\times S^2$. In~\cite{Bettiol17} he used Hoelzel's surgery
criterion~\cite{Hoelzel} to prove that the condition is preserved under
connected sums, and he classified, up to homeomorphism, the closed simply
connected four-manifolds that admit it; the list is the same as for
positive scalar curvature once the smooth structure is allowed to change.
Beyond the simply connected case, he asked the following
question~\cite[Section~3]{Bettiol17}.
\begin{question}
Is every finitely presented group the fundamental group of a closed
four-manifold with positive biorthogonal curvature?
\end{question}
\noindent Bettiol pointed out that an affirmative answer would follow if
the condition were stable under the corresponding surgeries, to which
Hoelzel's criterion does not apply. He singled out $S^2\times T^2$ as
the crucial test case, and the question whether $S^2\times T^2$ admits
positive biorthogonal curvature is also Question~4 in the problem list of
Morgan and Pansu~\cite{MP}. We settled this test case in~\cite{CZprod} by
constructing explicit metrics of positive biorthogonal curvature on
$S^2\times T^2$, and more generally on every smooth oriented
$S^2$-bundle over a closed oriented surface.  Since $\Z^2$ is not virtually free, the
biorthogonal and the isotropic conditions already part ways on
$S^2\times T^2$. The purpose of the present paper is to answer Bettiol's
question in full.

\subsection{Main results}

\begin{theorem}\label{thm:realization}
A group $\Gamma$ is isomorphic to the fundamental group of a closed
connected oriented smooth four-manifold admitting a metric of positive
biorthogonal curvature if and only if $\Gamma$ is finitely presented.
\end{theorem}

Only the ``if'' part has geometric content, since a closed smooth manifold
has a finite triangulation and its fundamental group is therefore finitely
presented. Thus positive biorthogonal curvature and positive scalar
curvature allow exactly the same fundamental groups.

Our construction also gives some control over the topology of the
manifold. For a finitely presented group $\Gamma$ we put
$b_1(\Gamma)=\dim_{\Q}H_1(\Gamma;\Q)$. For a closed oriented
four-manifold $M$ we write $\chi(M)$ for its Euler characteristic,
$b_j(M)$ for its rational Betti numbers, and $\sign(M)$ for the signature
of the intersection form on $H_2(M;\Q)$. Recall that an oriented manifold
is spin if its second Stiefel--Whitney class vanishes, and that a manifold
is \emph{stably parallelizable} if the sum of its tangent bundle and a
trivial line bundle is trivial. A stably parallelizable closed oriented
four-manifold is spin, and its signature vanishes because its first
Pontryagin class is zero~\cite{MilnorStasheff}.

\begin{corollary}\label{cor:topology}
Suppose that $\Gamma$ has a presentation with $p$ generators and $q$
relators, none of which is trivial in the free group on the generators.
Then the manifold $M$ in Theorem~\ref{thm:realization} can be chosen to
be stably parallelizable, with
\begin{align*}
 \chi(M)=2(1-p+q), \qquad
 b_2(M)=2\bigl(q-p+b_1(\Gamma)\bigr).
\end{align*}
\end{corollary}

The manifold $M$ is the boundary of a five-dimensional handlebody that is
built from the presentation, and it is of the kind that Carr used for
positive scalar curvature; see Remarks~\ref{rem:bounding}
and~\ref{rem:carr}. For positive scalar curvature the corresponding
statement was recorded by Mantione and Torres~\cite[Theorem~6]{MT}.

Behind Theorem~\ref{thm:realization} is a surgery result, which we state
separately. Let $\gamma$ be an embedded circle in an oriented
four-manifold $M$. Its normal bundle is trivial and has exactly two
oriented trivializations up to homotopy, the \emph{normal framings} of
$\gamma$.
\emph{Surgery along $\gamma$} removes an open tubular neighborhood
$S^1\times\operatorname{int}D^3$ of $\gamma$ and glues in $D^2\times S^2$
along the common boundary $S^1\times S^2$; the result depends on $\gamma$
and on the chosen framing.

\begin{theorem}\label{thm:surgery}
Let $M$ be a closed connected oriented smooth four-manifold that admits a
metric of positive biorthogonal curvature. Then every manifold obtained
from $M$ by surgery along finitely many pairwise disjoint embedded
circles, with arbitrary normal framings, also admits a metric of positive
biorthogonal curvature.
\end{theorem}

Bettiol's connected-sum theorem~\cite[Proposition~3.1]{Bettiol17} takes
care of surgery in codimension four. Together with
Theorem~\ref{thm:surgery} it shows that the class of closed oriented
four-manifolds that admit positive biorthogonal curvature is closed under
surgery of codimension at least three. The class of
manifolds that admit positive scalar curvature has the same property by
the theorem of Gromov--Lawson and Schoen--Yau. There is a difference,
however. The surgery theorem for positive scalar curvature modifies a
given metric near the surgery sphere only, whereas
Theorem~\ref{thm:surgery} is a statement about diffeomorphism types. In
its proof we first move the circle, by an isotopy of a stabilization of
$M$, to a closed geodesic of a suitable metric, and we do not claim that a
given metric can be modified near a given circle only; compare
Subsection~\ref{subsec:questions}.

Theorem~\ref{thm:surgery} does not follow from the general surgery theory
for curvature conditions. The model for a surgery of codimension three is
$S^2\times\R^2$, and its mixed planes are flat and have flat complements.
Its curvature operator therefore lies on the boundary of the cone of
operators with positive biorthogonal curvature, whereas Hoelzel's
criterion~\cite[Theorem~B]{Hoelzel} requires it to lie in the interior.
We get around this by twisting the handle; see
Subsection~\ref{subsec:outline}.

Attaching a two-handle to a five-manifold changes its boundary by surgery
along a circle. Theorem~\ref{thm:surgery} therefore has the following
consequence, which contains the existence statement of
Theorem~\ref{thm:realization}.

\begin{corollary}\label{cor:handlebody}
Let $W$ be a compact connected oriented smooth five-manifold that has a
handle decomposition with handles of index at most two. Then $\partial W$
admits a metric of positive biorthogonal curvature.
\end{corollary}

We mention some examples. Recall that the \emph{double} of a compact
manifold $X$ with boundary is the union of two copies of $X$ along their
common boundary, and that the \emph{spun manifold} of a closed connected
three-manifold $Y$ is the boundary of $Y_0\times D^2$, where $Y_0$ is the
complement of an open ball in $Y$ and $D^2$ is the closed two-dimensional
disk.

\begin{corollary}\label{cor:examples}
The following manifolds admit metrics of positive biorthogonal curvature.
\begin{enumerate}[label=\textup{(\alph*)}]
\item Every oriented $S^2$-bundle over a closed oriented surface.
\item The product $Y\times S^1$, for every closed orientable
three-manifold $Y$ that admits a metric of positive scalar curvature.
\item The double of every compact connected oriented four-manifold that
has a handle decomposition with handles of index at most two.
\item The spun manifold of every closed connected oriented three-manifold.
\end{enumerate}
\end{corollary}

Part~(a) is the main result of~\cite{CZprod}, and part~(b) with
$Y=S^2\times S^1$ gives a second proof for $S^2\times T^2$. Parts~(a),
(c) and~(d) are special cases of Corollary~\ref{cor:handlebody}.

If we allow connected sums with copies of $S^2\times S^2$, we can also
treat manifolds that are not given to us by a surgery description. Let
$N$ be a closed connected oriented four-manifold with fundamental group
$\Gamma$, and let $c$ be a classifying map $N\to B\Gamma$ of its
universal cover, where $B\Gamma$ denotes an Eilenberg--MacLane space of
type $K(\Gamma,1)$.

\begin{corollary}\label{cor:stable}
Suppose that $N$ satisfies one of the following conditions.
\begin{enumerate}[label=\textup{(\alph*)}]
\item The universal cover of $N$ is not spin, and
$c_*[N]=0$ in $H_4(\Gamma;\Z)$.
\item $N$ is spin, and for some spin structure the pair $(N,c)$ is
null-bordant in the spin bordism group $\Omega_4^{\mathrm{Spin}}(B\Gamma)$.
\end{enumerate}
Then the connected sum of $N$ with $r$ copies of $S^2\times S^2$ admits a
metric of positive biorthogonal curvature for all sufficiently large $r$.
\end{corollary}

The condition~(a) of Corollary~\ref{cor:stable} holds, for instance, for every $N$ with non-spin universal
cover whose fundamental group satisfies $H_4(\Gamma;\Z)=0$; this includes
free groups, surface groups and finite cyclic groups.

Theorem~\ref{thm:realization} does not say that every closed four-manifold
of positive scalar curvature admits positive biorthogonal curvature, and
we do not know whether this is true.

\begin{question}\label{q:psc}
Is there a closed four-manifold that admits a metric of positive scalar
curvature but no metric of positive biorthogonal curvature?
\end{question}

In the simply connected case, Bettiol's classification shows that such a
manifold cannot be detected by its homeomorphism type. Our results add
two restrictions. By Theorem~\ref{thm:surgery}, the two classes of
manifolds are closed under the same surgeries, so no example can be
produced, or excluded, by surgery alone. By
Corollary~\ref{cor:stable}, in the range covered there the two conditions
cannot be separated by any invariant that survives connected sums with
$S^2\times S^2$. An example answering Question~\ref{q:psc} affirmatively
would therefore have to rest on an obstruction to positive biorthogonal
curvature that is not an obstruction to positive scalar curvature and
that is unstable in this sense. We know of no such obstruction; see
Subsection~\ref{subsec:questions} for further discussion.

Finally, we study the minimal Euler characteristic $q_\perp(\Gamma)$ of a
closed connected oriented four-manifold with fundamental group $\Gamma$
that admits positive biorthogonal curvature. This is the analogue, for
our curvature condition, of the Hausmann--Weinberger invariant
$q(\Gamma)$~\cite{HW85}, the minimum over all closed oriented
four-manifolds with fundamental group $\Gamma$, and of the invariant
$q_{>0}(\Gamma)$ of Mantione and Torres~\cite{MT}, the minimum over those
that admit positive scalar curvature. We show that
\[
 q(\Gamma)\le q_{>0}(\Gamma)\le q_\perp(\Gamma)\le2-2\defic(\Gamma),
\]
where $\defic(\Gamma)$ is the deficiency of $\Gamma$
(Corollary~\ref{cor:deficiency}). These four numbers agree if there is a
field $\F$ such that $H_4(\Gamma;\F)=0$ and
$\defic(\Gamma)=b_1(\Gamma;\F)-b_2(\Gamma;\F)$, where
$b_j(\Gamma;\F)=\dim_{\F}H_j(\Gamma;\F)$
(Corollary~\ref{cor:minimal}). For $q_{>0}(\Gamma)$ an equality of this
kind is due to Mantione and Torres~\cite[Theorem~B]{MT}, and the new
information is that the minimum is attained by manifolds of positive
biorthogonal curvature. The corollary applies to the fundamental group of
a finite aspherical two-dimensional complex $X$, where the common value is
$2\chi(X)$ (Corollary~\ref{cor:aspherical-complex}); for $q(\Gamma)$ this
is a theorem of Hillman~\cite{Hillman94}. It also applies to the
fundamental group of a closed orientable aspherical three-manifold, where
the common value is $2$. For the fundamental group of a closed oriented
surface $\Sigma$ of positive genus we obtain $q_\perp=2\chi(\Sigma)$, and
the minimum is attained by both $S^2$-bundles over the surface. The
invariants $q$ and $q_\perp$ do not agree in general, since $q(\Z^4)=0$
and $q_\perp(\Z^4)=6$ (Remark~\ref{rem:Z4}). We do not know whether
$q_{>0}$ and $q_\perp$ can differ (Question~\ref{q:qperp}).

\subsection{Outline of the construction}\label{subsec:outline}

In dimension four, surgery along a circle replaces
$S^1\times\operatorname{int}D^3$ by $D^2\times S^2$, where $D^j$ denotes
the closed $j$-dimensional disk. As explained above, the standard model
$S^2\times\R^2$ is of no use to us. The problem is not an artifact of the
model. Let $r$ and $\theta$ be polar coordinates on the disk, and consider
an untwisted warped product $dr^2+f(r)^2d\theta^2+h(r)^2g_{S^2}$ on
$D^2\times S^2$. The planes spanned by $\partial_r$ and a direction
tangent to the sphere have biorthogonal curvature $-(fh')'/(2fh)$. The
function $fh'$ vanishes at the center of the disk and is positive on the
collar of a round tube. Hence $(fh')'$ is positive somewhere, and the
biorthogonal curvature of these planes is negative there.

We compensate for this by letting the sphere rotate as we go around the
disk. With spherical coordinates $\vartheta,\phi$ on $S^2$, we consider
\begin{equation}\label{eq:intro-ansatz}
 g=dr^2+f(r)^2d\theta^2+
 h(r)^2\bigl\{d\vartheta^2+\sin^2\vartheta\,(d\phi+A(r)\,d\theta)^2\bigr\},
\end{equation}
where the function $A$ is the coefficient of a rotational connection.
Its curvature enters through the quantity $\kappa=h^2A'/f$, and it
contributes the positive term $\kappa^2\sin^2\vartheta/(8h^2)$ to the
biorthogonal curvature of the critical planes. This term vanishes at the
poles of the sphere, which are fixed by the rotation. We then use the
fact that $k_g$ transforms under a conformal change of the metric exactly
like scalar curvature in dimension four. The conformal factor
$u=1+\kappa^2\cos^2\vartheta/24$ redistributes the positive term evenly
and turns it into $\kappa^2/(12h^2)$, uniformly on the sphere. If $h$
increases sufficiently slowly, this term dominates the negative
contribution of the convexity of $h$. Once the sphere is large enough, we
make $h$ concave, remove the connection curvature, and match the metric
exactly to the collar of a round tube. The handle obtained in this way is
very long and thin, but all its parameters are explicit and finite. The
precise statement is Theorem~\ref{surg:handle}.

To use the handle we have to bring the metric near the surgery circle
into the form of a round tube. We can do this when the circle is a closed
geodesic. The given metric and the round model then agree to first order
along the circle, and the local flexibility theorem of B\"ar and
Hanke~\cite{BarHanke} allows us to replace one by the other inside the
cone of positive biorthogonal curvature. The round tube remembers the
normal holonomy of the geodesic, an arbitrary rotation, and exact matching
has to take it into account; a portion of the handle of constant radius
adjusts the total rotation of the connection, and an integral twist of
the boundary realizes either normal framing.

It remains to produce the geodesics. For Theorem~\ref{thm:realization} we
replace each relator $r_j$ of a presentation by the pair $r_jy_j$, $y_j$,
where $y_j$ is a new generator. The new relators are primitive elements
of the enlarged free group, so that the corresponding free homotopy
classes on a connected sum of copies of $S^1\times S^3$ contain prime
closed geodesics; finite perturbation lemmas of
Rademacher~\cite{Rademacher} make these geodesics embedded and pairwise
disjoint. For Theorem~\ref{thm:surgery} we use the same device
geometrically. We take the connected sum of $M$ with two copies of
$S^1\times S^3$, slide the given circle over them, and use the fact that
homotopic circles in a four-manifold are isotopic.

\subsection*{Organization}
In Section~\ref{sec:preliminaries} we collect the curvature identities that we need.
In Section~\ref{sec:preparation} we prepare the surgery circles and their
tubular neighborhoods. Section~\ref{sec:handle} contains the construction of the
handle and is the technical core of the paper. The main results are proved
in Section~\ref{sec:realization}. Section~\ref{sec:discussion} contains the results
on Euler characteristics, the proof of Corollary~\ref{cor:stable}, and some
open questions.

\section{Preliminaries on biorthogonal curvature}\label{sec:preliminaries}

Let $(M^4,g)$ be an oriented Riemannian four-manifold with Levi-Civita
connection $\nabla$. We use the conventions
\begin{align*}
 R(X,Y)Z&=\nabla_X\nabla_YZ-\nabla_Y\nabla_XZ-\nabla_{[X,Y]}Z,\\
 \ip{X\wedge Y}{Z\wedge W}&=g(X,Z)g(Y,W)-g(X,W)g(Y,Z),
\end{align*}
and we define the curvature operator
$\Sc$ as the self-adjoint endomorphism of $\Lambda^2TM$ with
\begin{align*}
 \ip{\Sc(X\wedge Y)}{Z\wedge W}=g(R(X,Y)W,Z).
\end{align*}
We write $\Sc(\alpha,\beta)=\ip{\Sc\alpha}{\beta}$ for the associated
symmetric bilinear form, and we use $g$ to identify two-vectors with
two-forms, keeping the same notation for the curvature operator on
$\Lambda^2T^*M$. If $X$ and $Y$ are orthonormal, then
$K_g(\operatorname{span}\{X,Y\})=\Sc(X\wedge Y,X\wedge Y)$; in
particular, the round unit sphere has sectional curvature one. For an
orthonormal frame $e_1,\ldots,e_4$ we abbreviate
\[
 R_{ij,kl}=\Sc(e_i\wedge e_j,e_k\wedge e_l),
 \qquad\text{so that}\qquad R_{ij,ij}=K_g(e_i,e_j)
\]
is the sectional curvature of the plane spanned by $e_i$ and $e_j$.

\subsection{The Hodge decomposition}
The Hodge star $*$ is an isometric involution of $\Lambda^2T_pM$. Its
eigenspaces $\Lambda_p^\pm=\{\alpha\mid *\alpha=\pm\alpha\}$ are
three-dimensional and orthogonal to each other. If $e_1,\ldots,e_4$ is an
oriented orthonormal basis of $T_pM$, then
\begin{equation}\label{eq:hodge-basis}
\begin{aligned}
 \omega_1^\pm&=\tfrac1{\sqrt2}(e_1\wedge e_2\pm e_3\wedge e_4),\\
 \omega_2^\pm&=\tfrac1{\sqrt2}(e_1\wedge e_3\mp e_2\wedge e_4),\\
 \omega_3^\pm&=\tfrac1{\sqrt2}(e_1\wedge e_4\pm e_2\wedge e_3)
\end{aligned}
\end{equation}
are orthonormal bases of $\Lambda_p^\pm$. With respect to the
decomposition $\Lambda^2TM=\Lambda^+\oplus\Lambda^-$ the curvature
operator has the block form
\begin{align}\label{eq:curvature-blocks}
 \Sc=\begin{pmatrix}A_+&B\\ B^*&A_-\end{pmatrix},\qquad
 A_\pm=W_\pm+\frac{s_g}{12}\,I_{\Lambda^\pm},
\end{align}
where $W_\pm$ are the self-dual and the anti-self-dual Weyl curvature
operators and the off-diagonal block $B$ is given by the traceless Ricci
tensor; see Singer and Thorpe~\cite{ST69} and~\cite[Section~2]{CR14}. We
shall only need the diagonal blocks $A_\pm$, because of the following
formula, which reduces every curvature estimate in this paper to an
estimate for two symmetric $3\times3$ matrices. We write
$\lambda_{\min}$ for the smallest eigenvalue of a self-adjoint operator.

\begin{lemma}[\cite{CR14}]\label{lem:hodge-minimum}
For every oriented Riemannian four-manifold $(M,g)$ and every $p\in M$,
\begin{align}\label{eq:min}
 k_g(p)=\frac12\bigl(\lambda_{\min}(A_+(p))
                      +\lambda_{\min}(A_-(p))\bigr).
\end{align}
\end{lemma}

\begin{proof}
We recall the short argument. Let $\sigma$ be a two-plane with oriented
orthonormal basis $e_1,e_2$, complete it to an oriented orthonormal basis
of $T_pM$, and put $\alpha=e_1\wedge e_2$, so that $*\alpha=e_3\wedge e_4$
represents $\sigma^\perp$. Writing $\alpha=\alpha_++\alpha_-$ with
$\alpha_\pm=\frac12(\alpha\pm*\alpha)\in\Lambda^\pm$, we have
$*\alpha=\alpha_+-\alpha_-$, and therefore
\[
 K_g^\perp(\sigma)=\tfrac12\bigl(\Sc(\alpha,\alpha)+\Sc(*\alpha,*\alpha)\bigr)
 =A_+(\alpha_+,\alpha_+)+A_-(\alpha_-,\alpha_-),
 \qquad |\alpha_+|^2=|\alpha_-|^2=\tfrac12 .
\]
Conversely, every pair of vectors $\alpha_\pm\in\Lambda^\pm_p$ of norm
$1/\sqrt2$ arises in this way, since $\alpha_++\alpha_-$ is then a unit
two-vector with $\alpha\wedge\alpha=0$, hence decomposable. Minimizing
over the two factors separately gives \eqref{eq:min}.
\end{proof}

\subsection{Conformal changes}
For a function $u$ on $(M,g)$ we use the Laplacian
$\Delta_gu=\operatorname{tr}_g\nabla^2u$. Under a conformal change of the
metric, the Hessian terms in the transformation law of the sectional
curvature depend on the plane, but they cancel when one averages over a
plane and its orthogonal complement. As a consequence, $12\,k_g$
transforms exactly like the scalar curvature of a four-manifold. This is
the reason why the biorthogonal condition is accessible to conformal
methods, and it underlies the modified Yamabe invariants of
Costa~\cite{Costa12}.

\begin{lemma}[Conformal change]\label{lem:conformal}
Let $(M^4,g)$ be a Riemannian four-manifold and let $u$ be a positive
smooth function on $M$. For every two-plane $\sigma\subset T_pM$,
\begin{align*}
 K_{u^2g}^\perp(\sigma)
 =u^{-2}\left(K_g^\perp(\sigma)-\frac{\Delta_gu}{2u}\right).
\end{align*}
Consequently,
\begin{align}\label{eq:conformal}
 k_{u^2g}=u^{-2}\left(k_g-\frac{\Delta_gu}{2u}\right).
\end{align}
\end{lemma}

\begin{proof}
Write $u=e^{\varphi}$. For $g$-orthonormal vectors $X,Y$, the sectional
curvature of $e^{2\varphi}g$ is given by
\[
 e^{2\varphi}K_{e^{2\varphi}g}(X,Y)
 =K_g(X,Y)-\nabla^2\varphi(X,X)-\nabla^2\varphi(Y,Y)
  +(X\varphi)^2+(Y\varphi)^2-|\nabla\varphi|^2;
\]
see~\cite[Theorem~1.159]{B1987}. Orthogonal complements are the same for
the two metrics. Let $e_1,\ldots,e_4$ be a $g$-orthonormal basis with
$\sigma=\operatorname{span}\{e_1,e_2\}$. Adding the formulas for
$(e_1,e_2)$ and $(e_3,e_4)$, we obtain
\[
 e^{2\varphi}\bigl(K_{e^{2\varphi}g}(\sigma)+K_{e^{2\varphi}g}(\sigma^\perp)\bigr)
 =K_g(\sigma)+K_g(\sigma^\perp)-\Delta_g\varphi-|\nabla\varphi|^2,
\]
and $\Delta_g\varphi+|\nabla\varphi|^2=u^{-1}\Delta_gu$. See
also~\cite[Lemma~2.2]{CZprod}.
\end{proof}

\subsection{The cone of positive biorthogonal curvature}

In order to interpolate between two metrics we shall use the fact that
positive biorthogonal curvature is a convex condition on the curvature
tensor, a fact that also underlies Bettiol's connected-sum
theorem~\cite[Proposition~3.1]{Bettiol17}. We record it together with the
quantitative inequality behind it.

An \emph{algebraic curvature tensor} on a Euclidean vector space is a
tensor with the symmetries of a Riemannian curvature tensor, including
the first Bianchi identity. For such a tensor $T$ on a four-dimensional
Euclidean space we write $K_T(P)$ for its sectional curvature on a
two-plane $P$, and we put
\[
 k(T)=\min_{P}\frac{K_T(P)+K_T(P^\perp)}2,
\]
where $P$ runs through the Grassmannian of two-planes.

\begin{lemma}\label{lem:convex-cone}
The function $k$ is concave, that is,
\[
 k\bigl((1-a)T_0+aT_1\bigr)\ge (1-a)\,k(T_0)+a\,k(T_1)
 \qquad\text{for }0\le a\le1.
\]
The set $\{T\mid k(T)>0\}$ is an open, convex, $\mathrm O(4)$-invariant cone.
On a closed four-manifold, the set of metrics of positive biorthogonal
curvature is open in the $C^2$ topology.
\end{lemma}

\begin{proof}
For a fixed plane $P$, the quantity $(K_T(P)+K_T(P^\perp))/2$ is linear in
$T$, and a minimum of linear functions is concave. This gives the
inequality and the convexity of the cone. Since the Grassmannian is
compact, $k$ is continuous, indeed locally Lipschitz, and so the cone is
open; its invariance under $\mathrm O(4)$ and under multiplication by
positive constants is clear from the definition. Finally, the curvature
tensor of a metric depends continuously on the metric and its derivatives
up to order two, and on a closed manifold a metric with $k_g>0$ satisfies
$k_g\ge c$ for some constant $c>0$. This proves the last assertion.
\end{proof}

\section{Surgery circles and their tubular neighborhoods}
\label{sec:preparation}

In order to attach the handle of Section~\ref{sec:handle} we need the metric near
each surgery circle to be a round tube, in a sense made precise in
Lemma~\ref{lem:exact-tubes}. We achieve this for circles that are closed
geodesics. In this section we first show how to obtain closed geodesics
in the free homotopy classes that we need, and that they can be made
embedded and disjoint. We then replace the metric near such a geodesic
by a round tube, keeping track of the normal holonomy, and finally we
recall the effect of circle surgery on the fundamental group.

\subsection{Primitive relators and embedded geodesics}

For $n\ge1$ let $F_n$ be the free group on $n$ given generators, and let
$F_0$ be the trivial group. An element of a free group is
\emph{primitive} if it belongs to a free basis, and an element $w$ of a
group is a \emph{proper power} if $w=v^d$ for some element $v$ and some
integer $d\ge2$. A primitive element is nontrivial and is not a proper
power, because there is a homomorphism from the free group to $\Z$ that sends it
to $1$, whereas the image of $v^d$ is divisible by $d$.

Free homotopy classes of oriented loops in a connected manifold
correspond to conjugacy classes in its fundamental group, and reversing
the orientation of a loop replaces the conjugacy class by that of the
inverse. Thus a finite set of group elements that are pairwise
nonconjugate, even after inverting some of them, determines free homotopy
classes that are pairwise distinct up to reversal.

A length-minimizing closed geodesic in a class that is a proper power may
well run several times around a shorter geodesic, and two relators of a
presentation may be conjugate. The following elementary change of
presentation rules out both phenomena.

\begin{lemma}
\label{lem:primitive-presentation}
Every finitely presented group has a finite presentation whose relators
are primitive elements of the free group on the generators and are
pairwise nonconjugate, even after inversion.
\end{lemma}

\begin{proof}
Let
$\Gamma=\langle x_1,\ldots,x_p\mid r_1,\ldots,r_q\rangle$ be a finite
presentation. We freely reduce the words $r_j$ and discard those that
become empty, so that each $r_j$ is a nontrivial element of $F_p$. If
$q=0$ there is nothing to prove. Otherwise we introduce new generators
$y_1,\ldots,y_q$ and consider the presentation
\begin{equation}\label{eq:primitive-presentation}
 \left\langle
 x_1,\ldots,x_p,y_1,\ldots,y_q\;\middle|\;
 r_1y_1,\,y_1,\,\ldots,\,r_qy_q,\,y_q
 \right\rangle .
\end{equation}
It defines the same group $\Gamma$, since the relations $y_j=1$ eliminate the
new generators and turn $r_jy_j=1$ into $r_j=1$.
Each $y_j$ is primitive in $F_{p+q}$. The endomorphism of $F_{p+q}$ that
sends $y_j$ to $r_jy_j$ and fixes the other generators is an
automorphism, with inverse given by $y_j\mapsto r_j^{-1}y_j$; hence
$r_jy_j$ is primitive as well.

To compare the relators, consider the homomorphism $F_{p+q}\to\Z^q$ that
sends $y_j$ to the $j$-th standard basis vector $e_j$ and every $x_i$ to
zero. It maps both $r_jy_j$ and $y_j$ to $e_j$. Since conjugate elements
have the same image, relators with different indices are not conjugate,
and no relator is conjugate to the inverse of a relator. It remains to
show that $r_jy_j$ is not conjugate to $y_j$. The retraction
$F_{p+q}\to F_p$ that kills all $y_i$ maps these two elements to $r_j$
and to $1$, and $r_j$ is nontrivial in $F_p$.
\end{proof}

We next represent such classes by embedded closed geodesics. Recall that
a closed geodesic is \emph{prime} if it is not an iterate of a shorter
closed geodesic, and that two closed geodesics are \emph{geometrically
distinct} if their images differ.

\begin{lemma}
\label{lem:geodesic-representatives}
Let $(M,g)$ be a closed connected Riemannian manifold of dimension at
least three, and fix finitely many nontrivial free homotopy classes of
loops in $M$ that are not proper powers and are pairwise distinct up to
reversal. Then there is a smooth metric, arbitrarily close to $g$ in the
$C^2$ topology, for which the given classes are represented by pairwise
disjoint embedded closed geodesics. This metric can be chosen to agree
with $g$ outside finitely many balls of arbitrarily small radius. If $M$
has dimension four and $k_g>0$, the new metric can be chosen to have
positive biorthogonal curvature.
\end{lemma}

\begin{proof}
Since $M$ is compact, every nontrivial free homotopy class contains a
closed geodesic of minimal length. Let $c_1,\ldots,c_m$ be such
minimizers for the given classes. Each $c_i$ is prime, because a $d$-fold
iterate with $d\ge2$ represents a $d$-th power. The $c_i$ are also
geometrically distinct, because two prime closed geodesics with the same
image differ only by a shift or a reversal of the parameter, and they
would then represent classes that agree up to reversal.

Let $x$ be a double point of some $c_i$ or a common point of two different
$c_i$. The tangent lines at $x$ of the arcs of the $c_i$ through $x$ are
pairwise different, since two arcs with the same tangent line would
coincide by the uniqueness of geodesics, and this would contradict the
fact that the $c_i$ are prime and geometrically distinct. Hence such
points are isolated, and by compactness there are only finitely many of
them; see also~\cite[Lemma~2]{Rademacher}. For small $\eta>0$, the ball of
radius $2\eta$ around such a point $x$ meets the union of the $c_i$ in
finitely many geodesic arcs through $x$ only. As $\dim M\ge3$, we can
apply~\cite[Lemma~3]{Rademacher} to these arcs. In every neighborhood of
$g$ in the $C^2$ topology it provides a metric that agrees with $g$
outside a thin shell inside the ball, together with pairwise disjoint
geodesic arcs of the new metric that coincide with the original arcs near
the boundary of the ball. The new metric is smooth, because it is obtained
from $g$ by the construction of~\cite[Lemma~1]{Rademacher}, which uses
only smooth diffeomorphisms and smooth cutoff functions.

We choose the balls around the various points $x$ to be pairwise
disjoint and carry out this modification in each of them. Together with
the unchanged parts of the $c_i$, the new arcs form closed geodesics of
the new metric. These are embedded and pairwise disjoint, and each of them
is freely homotopic to the geodesic $c_i$ that it replaces, because the
two differ only inside balls.
Finally, if $\dim M=4$ and $k_g>0$, then a sufficiently small perturbation
in the $C^2$ topology preserves positive biorthogonal curvature by
Lemma~\ref{lem:convex-cone}.
\end{proof}

\begin{remark}\label{rem:finite-geodesic-perturbation}
For a $C^r$-generic metric with $r\ge2$, all prime closed geodesics are
simple and geometrically distinct closed geodesics are
disjoint~\cite[Theorem~1]{Rademacher}. We only need the finite version
above, in which the perturbation stays close to a given smooth metric of
positive biorthogonal curvature.
\end{remark}

\subsection{Replacing a metric by one with the same first jet}

We say that two metrics have the \emph{same first jet along a
submanifold} $N$ if, in local coordinates, their coefficients and the
first partial derivatives of their coefficients agree at every point of
$N$. This does not depend on the choice of coordinates.

\begin{lemma}
\label{lem:matched-jets}
Let $N$ be a compact embedded submanifold without boundary of a smooth
four-manifold $M$, and let $U\subset M$ be an open neighborhood of $N$.
Let $g_0$ and $g_1$ be smooth Riemannian metrics on $U$ that have the
same first jet along $N$, and suppose that $k_{g_0}>0$ on $U$ and
$k_{g_1}>0$ at the points of $N$. Then for every open neighborhood
$V\subset U$ of $N$ there is a smooth path $g(s)$, $0\le s\le1$, of
Riemannian metrics on $U$ with the following properties.
\begin{enumerate}[label=\textup{(\roman*)}]
\item $g(0)=g_0$, and $g(s)=g_0$ outside $V$ for all $s$.
\item $g(1)=g_1$ on a neighborhood of $N$.
\item Every $g(s)$ has positive biorthogonal curvature and has the same
first jet as $g_0$ along $N$.
\end{enumerate}
In particular, if $g_0$ is defined on all of $M$ and $k_{g_0}>0$,
the path extends to $M$ by $g_0$ outside $U$.
\end{lemma}

\begin{proof}
For $s\in[0,1]$ put $g_s=(1-s)g_0+sg_1$. In local coordinates, the
components of the curvature tensor $\operatorname{Rm}_g$ of type $(0,4)$
are the sum of a linear expression in the second derivatives of the
metric, with constant coefficients, and of a term that depends only on
the metric and its first derivatives. Since all the metrics $g_s$ have
the same first jet along $N$, we obtain
\begin{equation}\label{eq:matched-jets-affine-curvature}
 \operatorname{Rm}_{g_s}|_N
 =(1-s)\operatorname{Rm}_{g_0}|_N
       +s\operatorname{Rm}_{g_1}|_N .
\end{equation}
The metrics $g_s$ induce the same inner product on $TM|_N$, so
Lemma~\ref{lem:convex-cone} applies and shows that $k_{g_s}>0$ at the
points of $N$ for every $s$. Since $N\times[0,1]$ is compact, there is an
open neighborhood $W$ of $N$ with compact closure contained in $V$ such
that $k_{g_s}>0$ on $W$ for every $s$.

We now apply the local flexibility theorem of B\"ar and
Hanke~\cite[Theorem~1.2]{BarHanke}. In their notation, the ambient
manifold is $U$, the closed subset is $N$, and its open neighborhood is
$W$. The fiber bundle is the bundle $\operatorname{Sym}^2T^*U$ of
symmetric bilinear forms, and the partial differential relation is the
subset of its second jet bundle that consists of the jets of positive
definite forms with positive biorthogonal curvature; it is open, because
the curvature depends continuously on the two-jet of the metric. The
section $g_0$ solves this relation on $U$, the path $s\mapsto g_s|_W$
consists of solutions on $W$ and starts at $g_0|_W$, and the one-jets of
$g_s$ along $N$ do not depend on $s$. Since $g_0$ is smooth and the path
$s\mapsto g_s|_W$ is smooth, the regularity statement that is part
of~\cite[Theorem~1.2]{BarHanke} applies. The theorem therefore gives a
smooth path $g(s)$ of solutions on $U$ with $g(0)=g_0$ that agrees with
$g_s$ on a neighborhood of $N$ independent of $s$ and with $g_0$ outside
$W$. This path has all the required properties.
\end{proof}

\subsection{Round tubes and normal holonomy}

We now apply Lemma~\ref{lem:matched-jets} to an embedded closed geodesic.
The model metric is the product of a line with a round three-ball,
written in Fermi coordinates with respect to a \emph{parallel} normal
frame. Such a frame need not close up after one turn around the geodesic;
the resulting rotation is the normal holonomy, and it becomes part of the
model.

\begin{lemma}
\label{lem:exact-tubes}
Let $(M^4,g)$ be an oriented Riemannian four-manifold with $k_g>0$, let
$\gamma$ be an embedded closed geodesic of length $L$, and let $K>0$.
Then there are a number $\rho_0>0$ with $\sqrt K\rho_0<\pi/2$ and a smooth
path of metrics of positive biorthogonal curvature that starts at $g$, is
constant outside an arbitrarily small neighborhood of $\gamma$, and ends
at a metric that near $\gamma$ is isometric to
\begin{equation}\label{eq:parallel-round-tube}
 dt^2+d\rho^2+s_K(\rho)^2g_{S^2},\qquad
 s_K(\rho)=\frac{\sin(\sqrt K\,\rho)}{\sqrt K},
\end{equation}
on the quotient of $\R\times B^3_{\rho_0}$ by the identification
$(t+L,x)\sim(t,Hx)$. Here $B^3_{\rho_0}\subset\R^3$ is the ball of radius
$\rho_0$ around the origin, $\rho=|x|$, the metric $g_{S^2}$ is the round
metric of the unit sphere, and $H\in\mathrm{SO}(3)$ is the holonomy of the
normal bundle of $\gamma$ with respect to $g$. All metrics of the path have the same first
jet along $\gamma$; in particular, $\gamma$ remains a geodesic of length
$L$ with normal holonomy $H$. Finitely many pairwise disjoint embedded
closed geodesics can be treated simultaneously.
\end{lemma}

\begin{proof}
\step{Step 1. Fermi coordinates.}
We parametrize $\gamma$ by arc length on $\R$, so that $\gamma$ is
$L$-periodic, and we choose an oriented orthonormal frame
$E_1(t),E_2(t),E_3(t)$ of the normal bundle that is parallel along
$\gamma$. Then $E(t+L)=E(t)H$ for a matrix $H\in\mathrm{SO}(3)$, where
$E(t)$ stands for the row $(E_1(t),E_2(t),E_3(t))$. For
$x=(x^1,x^2,x^3)\in\R^3$ small we put
\[
 \Phi(t,x)=\exp_{\gamma(t)}\Bigl(\sum_{a=1}^3x^aE_a(t)\Bigr).
\]
Then $\Phi(t+L,x)=\Phi(t,Hx)$, and since $\gamma$ is compact and
embedded, $\Phi$ induces a diffeomorphism from the quotient of
$\R\times B^3_{\rho_0}$ by $(t+L,x)\sim(t,Hx)$ onto a tubular
neighborhood of $\gamma$, provided $\rho_0$ is small.

\step{Step 2. The first jet of $g$.}
We claim that in the coordinates $(t,x)$ the coefficients of $g$ satisfy
\begin{equation}\label{eq:fermi-first-jet}
 g_{ij}(t,0)=\delta_{ij},\qquad
 \partial_k g_{ij}(t,0)=0
\end{equation}
for all indices. Along $x=0$ the coordinate vector fields are
$\dot\gamma,E_1,E_2,E_3$, which are orthonormal. Their covariant
derivatives in the direction of $\partial_t$ vanish there, because
$\gamma$ is a geodesic and the $E_a$ are parallel; since coordinate vector
fields commute, we also have
$\nabla_{\partial_a}\partial_t=\nabla_{\partial_t}\partial_a=0$ along
$x=0$. Finally, for fixed $t$ and $v\in\R^3$ the curve
$s\mapsto\Phi(t,sv)$ is a geodesic, which gives
$\sum_{a,b}v^av^b\,\nabla_{\partial_a}\partial_b=0$ at $(t,0)$ for all
$v$ and hence $\nabla_{\partial_a}\partial_b=0$ there. Thus all
Christoffel symbols vanish along $\gamma$, and \eqref{eq:fermi-first-jet}
follows from the compatibility of $\nabla$ with the metric.

\step{Step 3. The model.}
Let $g_1$ be the metric \eqref{eq:parallel-round-tube}. It is the product
of the line with a geodesic ball in the three-sphere of curvature $K$,
written in normal coordinates, and it descends to the quotient because the
round ball is invariant under $\mathrm{SO}(3)$. Since
$s_K(\rho)=\rho-\frac K6\rho^3+O(\rho^5)$, the metric $g_1$ is smooth
along the core $x=0$ and also satisfies \eqref{eq:fermi-first-jet}. Hence
$g$ and $g_1$ have the same first jet along $\gamma$.

The biorthogonal curvature of $g_1$ is constant. Indeed, let
$e_0=\partial_t$. For a two-plane $P$ we have
$K_{g_1}(P)=K\bigl(1-|\operatorname{proj}_Pe_0|^2\bigr)$, where
$\operatorname{proj}_P$ is the orthogonal projection onto $P$; this is
seen by choosing an orthonormal basis of $P$ whose first vector is
orthogonal to $e_0$. Since the squared lengths of the projections of
$e_0$ to $P$ and to $P^\perp$ add up to one, we get
\[
 K_{g_1}(P)+K_{g_1}(P^\perp)=K
 \qquad\text{and hence}\qquad k_{g_1}=\frac K2>0 .
\]

\step{Step 4. Conclusion.}
We apply Lemma~\ref{lem:matched-jets} to $g$ and $g_1$, with $N=\gamma$
and with $V$ a small tube around $\gamma$. The final metric of the
resulting path agrees with $g_1$ on a neighborhood of $\gamma$, and after
decreasing $\rho_0$ we may assume that this neighborhood contains the tube
of radius $\rho_0$ and that $\sqrt K\rho_0<\pi/2$. The path has all the
properties stated. Indeed, since the metric and its first derivatives
along $\gamma$ do not change, neither do the Levi-Civita connection along
$\gamma$, the length of $\gamma$, and the parallel transport in its
normal bundle. For finitely many disjoint geodesics we use disjoint tubes.
\end{proof}

\subsection{Collar coordinates and normal framings}
We now rewrite the model \eqref{eq:parallel-round-tube} in coordinates
that are periodic along $\gamma$. This fixes the form of the collar to
which the handle has to be matched, and our sign conventions for the
holonomy.

We identify $\R^3$ with $\R\oplus\C$ and write the unit sphere as
\[
 S^2=\{(z,w)\in\R\oplus\C\mid z^2+|w|^2=1\}.
\]
After a constant rotation of the parallel frame $E$, the holonomy is a
rotation about the $z$-axis, and we can choose $\alpha\in\R$ with
\begin{equation}\label{eq:normal-holonomy-convention}
 H(z,w)=(z,e^{-2\pi i\alpha}w).
\end{equation}
The number $\alpha$ is determined up to an integer. Put $F=L/(2\pi)$ and
$\theta=t/F$, so that $\theta$ has period $2\pi$, and consider the map
$\Psi$ from $S^1\times(0,\rho_0)\times S^2$ to the quotient
$\bigl(\R\times B^3_{\rho_0}\bigr)/\!\sim$ that is given by
\[
 \Psi(\theta,\rho,(z,w))=\bigl[\,t,\ \rho\,(z,e^{i\alpha\theta}w)\,\bigr].
\]
It is well defined, because increasing $\theta$ by $2\pi$ is compensated
by the identification with $H$, and it is a diffeomorphism onto the
complement of the core. In these coordinates the metric
\eqref{eq:parallel-round-tube} becomes
\begin{equation}\label{eq:periodic-round-tube}
 F^2d\theta^2+d\rho^2+
 s_K(\rho)^2\bigl(dz^2+|dw+i\alpha w\,d\theta|^2\bigr).
\end{equation}
In terms of the colatitude $\vartheta\in[0,\pi]$ and the longitude
$\phi\in\R/2\pi\Z$, defined by $z=\cos\vartheta$ and
$w=\sin\vartheta\,e^{i\phi}$, the last bracket equals
$d\vartheta^2+\sin^2\vartheta\,(d\phi+\alpha\,d\theta)^2$. We finally use
$h=s_K(\rho)$ as a radial coordinate. Since $\sqrt K\rho_0<\pi/2$, we have
$dh/d\rho=\cos(\sqrt K\rho)>0$ on $(0,\rho_0)$, and we arrive at the
following normal form.

\begin{lemma}\label{lem:collar}
In the situation of Lemma~\ref{lem:exact-tubes}, let $2\pi F$ be the
length of $\gamma$ and choose $\alpha$ as in
\eqref{eq:normal-holonomy-convention}. Then, with $h_1=s_K(\rho_0)$, the
complement of $\gamma$ in a tubular neighborhood has coordinates
$(\theta,h,(z,w))\in S^1\times(0,h_1)\times S^2$ in which the deformed
metric is
\begin{equation}\label{eq:round-tube-boundary-form}
 F^2d\theta^2+\frac{dh^2}{1-Kh^2}
 +h^2\bigl\{d\vartheta^2+
 \sin^2\vartheta\,(d\phi+\alpha\,d\theta)^2\bigr\}.
\end{equation}
The coordinate $h$ tends to zero at $\gamma$.
\end{lemma}

The connection term $\alpha\,d\theta$ can be removed locally by rotating
the sphere, but not globally unless $\alpha$ is an integer, because its holonomy
around $\gamma$ is the normal holonomy of the geodesic.

The coordinates of Lemma~\ref{lem:collar} also determine a normal framing
of $\gamma$, namely the homotopy class of the trivialization of the
normal bundle by the periodic frame $E(t)$ rotated by the angle
$\alpha\theta$ about the first axis. The choice of $\alpha$ within its
class modulo integers matters here. For an integer $n$, the diffeomorphism
\begin{equation}\label{eq:normal-gauge-twist}
 T_n(\theta,z,w)=(\theta,z,e^{in\theta}w)
\end{equation}
of $S^1\times S^2$ pulls the metric \eqref{eq:round-tube-boundary-form}
with coefficient $\alpha$ back to the same expression with coefficient
$\alpha+n$, since
$T_n^*(dw+i\alpha w\,d\theta)=e^{in\theta}\bigl(dw+i(\alpha+n)w\,d\theta\bigr)$.
The loop $\theta\mapsto$ (rotation by $n\theta$) represents the class of
$n$ modulo two in $\pi_1(\mathrm{SO}(3))\cong\Z/2$. Hence replacing
$\alpha$ by $\alpha+n$ does not change the holonomy, and it changes the
framing if and only if $n$ is odd. Normal holonomy and normal framing are
thus independent data. The former is the class of $\alpha$ modulo
integers, and the latter is a parity.

\begin{remark}
\label{rem:exact-collar-gluing}
Let $0<R<h_1$ and remove the open tube $\{h<R\}$ around $\gamma$. Suppose
that we are given a metric on $D^2\times S^2$ whose boundary collar is
parametrized by $R_-<h\le R$ and has the form
\eqref{eq:round-tube-boundary-form} with coefficient $\alpha+n$ instead of
$\alpha$. The map $T_n$ identifies this collar isometrically with the
region $\{R_-<h\le R\}$ of the tube. We can therefore glue $D^2\times S^2$
to the complement of the tube along $\{h=R\}$, and the resulting metric is
smooth, because on both sides of the seam it is given by the same formula in the
coordinate $h$, which increases from the handle into the exterior. The
resulting manifold is the result of surgery along $\gamma$, with the
framing determined by the coordinates if $n$ is even and with the other
framing if $n$ is odd.
\end{remark}

\subsection{The effect of circle surgery on the fundamental group}

For elements $a_1,\ldots,a_m$ of a group we write
$\langle\!\langle a_1,\ldots,a_m\rangle\!\rangle$ for the normal subgroup
that they generate. It depends only on the conjugacy classes of the
$a_j$, so the notation makes sense for free homotopy classes of loops.

\begin{lemma}
\label{lem:circle-surgery}
Let $M$ be a closed connected oriented smooth four-manifold, and let
$\gamma_1,\ldots,\gamma_m$ be pairwise disjoint embedded circles in $M$.
For any choice of normal framings, the manifold $M'$ obtained from $M$ by
surgery along $\gamma_1,\ldots,\gamma_m$ is a closed connected oriented
smooth four-manifold with
\begin{equation}\label{eq:circle-surgery-fundamental-group}
 \pi_1(M')\cong
 \pi_1(M)/\langle\!\langle
 [\gamma_1],\ldots,[\gamma_m]
 \rangle\!\rangle .
\end{equation}
More precisely, if $W$ denotes the trace of the surgeries, a compact
oriented five-manifold with boundary $M'\sqcup -M$, then the inclusion
$M'\subset W$ induces an isomorphism of fundamental groups and the
inclusion $M\subset W$ induces the quotient map in
\eqref{eq:circle-surgery-fundamental-group}.
\end{lemma}

\begin{proof}
The trace $W$ is obtained from $M\times[0,1]$ by attaching one
two-handle $D^2\times D^3$ along a tubular neighborhood of each
$\gamma_j\times\{1\}$. Seen from $M'$, it is obtained from
$M'\times[0,1]$ by attaching handles of index three. Up to homotopy,
attaching a two-handle amounts to attaching a two-cell along the
corresponding circle, and attaching a three-handle to attaching a
three-cell. The assertions therefore follow
from~\cite[Proposition~1.26]{Hatcher}.
\end{proof}

\section{A handle with positive biorthogonal curvature}\label{sec:handle}

In this section we construct the metric on $D^2\times S^2$ that replaces
a round tube. Its boundary collar is the collar
\eqref{eq:round-tube-boundary-form} of a closed geodesic of length
$2\pi F$ whose normal holonomy is given by an arbitrary real number
$\alpha$. The integer $n$ in the statement is produced by the
construction, but we can prescribe its parity, and hence the framing.

\begin{theorem}\label{surg:handle}
Let $F,K>0$, let $R>0$ satisfy $KR^2\le 1/100$, let $\alpha\in\R$, and
let $\nu\in\{0,1\}$. Then there exist an integer $n\equiv\nu\pmod 2$, a
number $R_-\in(0,R)$, and a smooth Riemannian metric $\widehat g$ on
$D^2\times S^2$ with $k_{\widehat g}>0$ whose boundary has a collar
parametrized by $(\theta,h,(z,w))\in S^1\times(R_-,R]\times S^2$ on which
\begin{equation}\label{surg:outer-model}
 \widehat g=F^2d\theta^2+\frac{dh^2}{1-Kh^2}
 +h^2\bigl\{d\vartheta^2+
       \sin^2\vartheta\,(d\phi+(\alpha+n)\,d\theta)^2\bigr\}.
\end{equation}
Here $\theta$ has period $2\pi$, the boundary is $\{h=R\}$, and
$(\vartheta,\phi)$ are the spherical coordinates on
$S^2=\{(z,w)\in\R\oplus\C\mid z^2+|w|^2=1\}$ given by $z=\cos\vartheta$ and
$w=\sin\vartheta\,e^{i\phi}$.
\end{theorem}

By Remark~\ref{rem:exact-collar-gluing}, the metric $\widehat g$ can be
glued to the complement of the tube $\{h<R\}$ in
Lemma~\ref{lem:collar} by means of the map $T_n$ of
\eqref{eq:normal-gauge-twist}, and the two values of $\nu$ give the two
normal framings. The metric that we construct is invariant under the
rotations of the disk and under the rotations of the sphere about the
$z$-axis, and the integer $n$ depends only on $F$, $K$, $R$, $\alpha$,
$\nu$ and on the choice of a cutoff function. The rest of this section
contains the proof of Theorem~\ref{surg:handle}. We begin with an
overview.

Let $r$ and $\theta$ be polar coordinates on the disk. The simplest
candidate for a metric on $D^2\times S^2$ is the warped product
\[
 dr^2+f(r)^2\,d\theta^2+h(r)^2g_{S^2},
\]
in which $f$ is the radius of the $\theta$-circle and $h$ is the radius
of the sphere. Smoothness at the center of the disk forces $f(0)=0$,
$f'(0)=1$ and $h'(0)=0$. On the collar \eqref{surg:outer-model}, on the
other hand, $f=F$ and $h'^2=1-Kh^2$ is close to one. For the warped
product this is fatal. We shall see in \eqref{surg:untwisted} that its
biorthogonal curvature satisfies $k\le-(fh')'/(2fh)$. Since $fh'$ vanishes
at $r=0$ and is positive on the collar, the right-hand side is negative
somewhere. The planes responsible for this are the ones spanned by the
radial direction and a direction tangent to the sphere, together with
their complements, and they are the planes that are flat in the model
$S^2\times\R^2$. In our construction $f$ will be constant outside a
neighborhood of the center, and the relevant quantity is the
\emph{defect}
\[
 \delta=1-h'^2 .
\]
It has to decrease from $1$ to a small number, and so $h$ has to be convex
on some interval. Our construction remedies the resulting loss of
positivity by means of three devices.

\smallskip\noindent
\emph{A rotational connection.} We let the sphere rotate about its axis
as $\theta$ varies, and we consider the metric \eqref{surg:ansatz} below.
The function $\kappa=h^2A'/f$ is a normalized curvature of this
connection. By Lemma~\ref{lem:hodge-blocks}, it adds the term
$\kappa^2(1-z^2)/(8h^2)$ to the biorthogonal curvature, where
$z=\cos\vartheta$ is the height function on the sphere. This term is
positive except at the poles $z=\pm1$, which are the fixed points of the
rotation.

\smallskip\noindent
\emph{A conformal change.} Since $k_g$ transforms under conformal changes
like a scalar curvature (Lemma~\ref{lem:conformal}), we can use a
conformal factor to move positivity from the equator to the poles. The
factor $u=1+\kappa^2z^2/24$ is chosen in such a way that the dependence
on $z$ cancels in the leading terms. After the conformal change, the
connection contributes $\kappa^2/(12h^2)$ uniformly on the sphere
(Lemma~\ref{lem:conformal-identity}).

\smallskip\noindent
\emph{Slow growth.} The convexity of $h$ costs the term $-u\,h''/(2h)$,
and the gain is $\kappa^2/(12h^2)$. We cannot take $\kappa$ large,
because the eigenvalue structure that we use requires, roughly,
$3\kappa<\delta$. We take $\kappa=\delta/6$, and positivity then holds as
long as $hh''\le\epsilon\delta^2/2$ for a small constant $\epsilon$. In
terms of the logarithmic radius $t=\log h$ this means
$d\delta/dt=-\epsilon\delta^2$, so that $1/\delta$ grows only linearly in
$t$, and the sphere has to grow by a factor of about
$\exp(1/(\epsilon\delta))$ before the defect has come down to $\delta$.
This is why the neck of our handle is very long. It is finite, however,
and no limiting argument is needed.

\smallskip
Once the defect is small, we let it increase again. Now $h$ is concave,
which is favorable, and we can switch off the connection and steer the
defect to the value $Kh^2$ of the round collar. The total rotation
$\int A'\,dr$ that has accumulated along the way is then corrected by
inserting a cylinder of constant radius near the center, on which the
connection has constant curvature; its length is chosen so that the
boundary coefficient becomes $\alpha+n$ with $n$ of the prescribed parity.

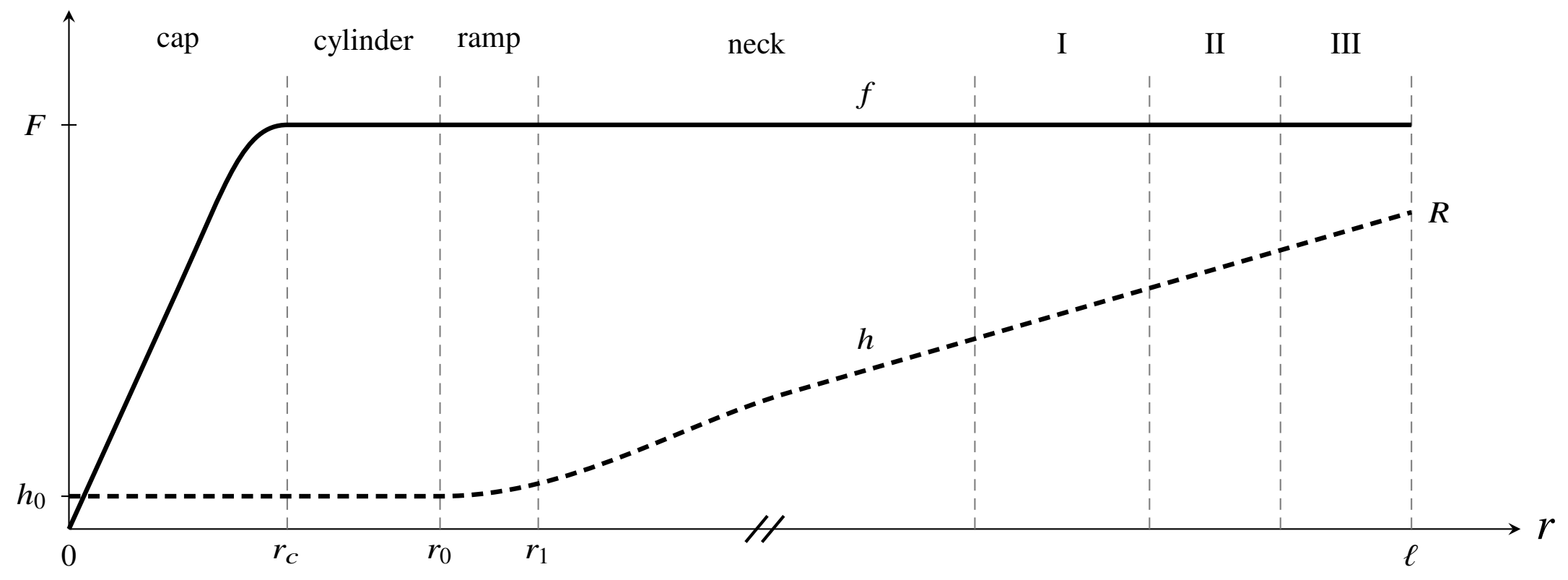
\begin{figure}[ht]
\centering
\begin{tikzpicture}[>=stealth,x=1cm,y=1cm]
 \draw[->] (0,0) -- (13.3,0) node[right] {$r$};
 \draw[->] (0,0) -- (0,4.75);
 \foreach \x in {2.0,3.4,4.3,8.3,9.9,11.1,12.3}
   \draw[densely dashed,gray] (\x,0) -- (\x,4.2);
 \node[below] at (0,0) {\scriptsize$0$};
 \node[below] at (2.0,0) {\scriptsize$r_c$};
 \node[below] at (3.4,0) {\scriptsize$r_0$};
 \node[below] at (4.3,0) {\scriptsize$r_1$};
 \node[below] at (12.3,0) {\scriptsize$\ell$};
 \draw[thick] (6.2,-0.12) -- (6.4,0.12);
 \draw[thick] (6.35,-0.12) -- (6.55,0.12);
 \draw[very thick] (0,0) -- (1.0,2.2) .. controls (1.45,3.19) and (1.6,3.7) .. (2.0,3.7) -- (12.3,3.7);
 \node[above] at (7.3,3.7) {\scriptsize$f$};
 \draw (-0.07,3.7) -- (0.07,3.7);
 \node[left] at (-0.05,3.7) {\scriptsize$F$};
 \draw[very thick,densely dashed] (0,0.3) -- (3.4,0.3)
   .. controls (4.6,0.3) and (5.6,0.96) .. (6.6,1.25) -- (12.3,2.9);
 \draw (-0.07,0.3) -- (0.07,0.3);
 \node[left] at (-0.05,0.3) {\scriptsize$h_0$};
 \node[right] at (12.3,2.9) {\scriptsize$R$};
 \node[above] at (7.3,1.5) {\scriptsize$h$};
 \node at (1.0,4.45) {\scriptsize cap};
 \node at (2.7,4.45) {\scriptsize cylinder};
 \node at (3.85,4.45) {\scriptsize ramp};
 \node at (6.3,4.45) {\scriptsize neck};
 \node at (9.1,4.45) {\scriptsize I};
 \node at (10.5,4.45) {\scriptsize II};
 \node at (11.7,4.45) {\scriptsize III};
\end{tikzpicture}
\caption{The radial profiles of the handle (schematic and not to scale).
The solid curve is the radius $f$ of the $\theta$-circle, the dashed
curve is the radius $h$ of the sphere. The disk closes up on the cap; the
cylinder has adjustable length and is used to tune the holonomy; the
sphere then grows from $h_0$ to $R$ along the ramp, the neck and
Stages~I--III. In reality $h_0$ is smaller than $R$ by many orders of
magnitude, and $h'$ is very close to one from the middle of the neck on.}
\label{fig:handle-r}
\end{figure}

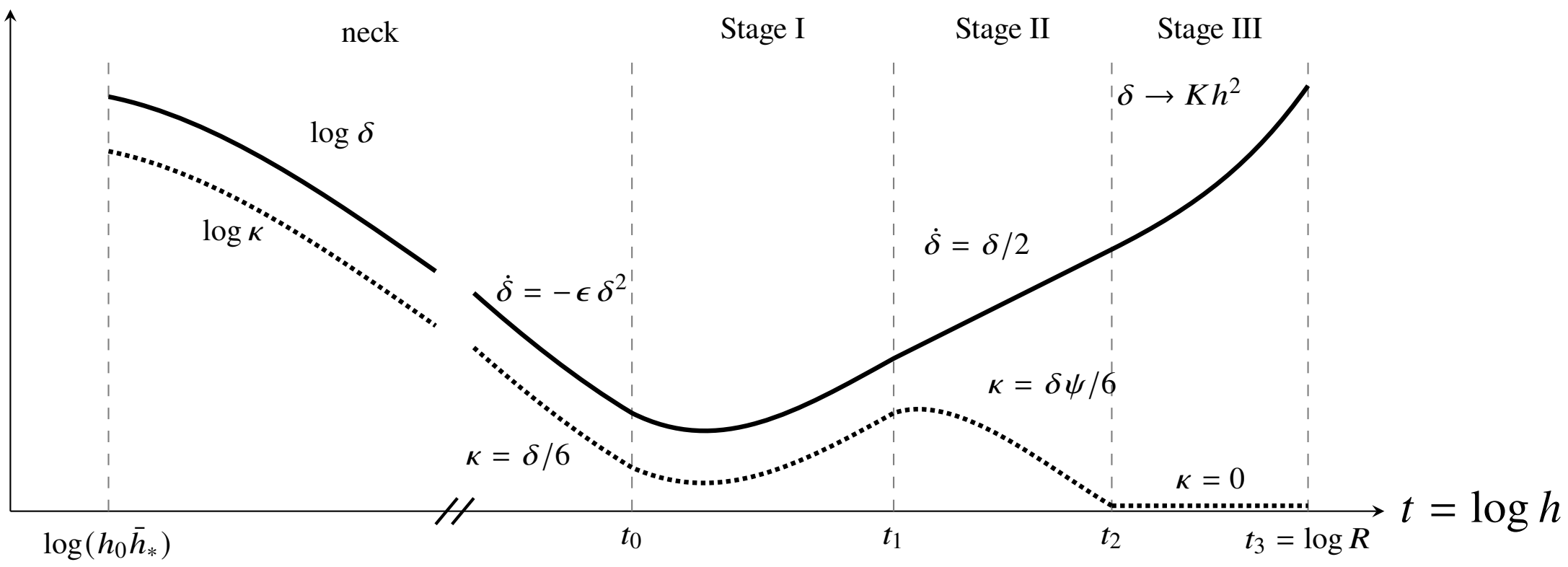
\begin{figure}[ht]
\centering
\begin{tikzpicture}[scale=1.0,>=stealth]
 \draw[->] (-0.3,0) -- (12.3,0) node[right] {$t=\log h$};
 \draw[->] (-0.3,0) -- (-0.3,4.6);
 \foreach \x in {0.6, 5.4, 7.8, 9.8, 11.6}
   \draw[dashed,gray] (\x,0) -- (\x,4.15);
 \node[below] at (0.6,0) {\scriptsize$\log(h_0\hn_*)$};
 \node[below] at (5.4,0) {\scriptsize$t_0$};
 \node[below] at (7.8,0) {\scriptsize$t_1$};
 \node[below] at (9.8,0) {\scriptsize$t_2$};
 \node[below] at (11.6,0) {\scriptsize$t_3=\log R$};
 \node at (3.0,4.4) {\scriptsize neck};
 \node at (6.6,4.4) {\scriptsize Stage I};
 \node at (8.8,4.4) {\scriptsize Stage II};
 \node at (10.7,4.4) {\scriptsize Stage III};
 \draw[thick] (3.6,-0.12) -- (3.8,0.12);
 \draw[thick] (3.75,-0.12) -- (3.95,0.12);
 \draw[very thick]
   (0.6,3.8) .. controls (1.6,3.6) and (2.6,2.9) .. (3.6,2.2);
 \draw[very thick]
   (3.95,2.0) .. controls (4.4,1.6) and (4.9,1.2) .. (5.4,0.9)
   .. controls (6.2,0.5) and (6.9,0.9) .. (7.8,1.4)
   -- (9.8,2.4)
   .. controls (10.6,2.8) and (11.1,3.2) .. (11.6,3.9);
 \node at (2.75,3.45) {\scriptsize$\log\delta$};
 \node at (4.75,2.05) {\scriptsize$\dot\delta=-\epsilon\delta^2$};
 \node at (8.55,2.45) {\scriptsize$\dot\delta=\delta/2$};
 \node[left] at (11.15,3.85) {\scriptsize$\delta\to Kh^2$};
 \draw[very thick,densely dotted]
   (0.6,3.3) .. controls (1.6,3.1) and (2.6,2.4) .. (3.6,1.7);
 \draw[very thick,densely dotted]
   (3.95,1.5) .. controls (4.4,1.1) and (4.9,0.7) .. (5.4,0.4)
   .. controls (6.2,0.05) and (6.9,0.4) .. (7.8,0.9)
   .. controls (8.4,1.1) and (9.2,0.4) .. (9.8,0.05);
 \draw[very thick,densely dotted] (9.8,0.05) -- (11.6,0.05);
 \node at (1.75,2.55) {\scriptsize$\log\kappa$};
 \node at (4.35,0.5) {\scriptsize$\kappa=\delta/6$};
 \node at (9.25,1.15) {\scriptsize$\kappa=\delta\psi/6$};
 \node at (10.7,0.3) {\scriptsize$\kappa=0$};
\end{tikzpicture}
\caption{The defect $\delta=1-h'^2$ (solid) and the normalized connection
curvature $\kappa$ (dotted) as functions of the logarithmic radius
$t=\log h$, on a logarithmic scale (schematic; the marks on the axis
indicate that most of the neck has been omitted). The cap, the cylinder
and the ramp, where $h$ is equal or very close to $h_0$, lie to the left
of the picture.}
\label{fig:handle-profiles}
\end{figure}

The pieces of the handle are listed in Table~\ref{tab:pieces} and
sketched in Figures~\ref{fig:handle-r} and~\ref{fig:handle-profiles}.

\begin{table}[ht]
\centering
\renewcommand{\arraystretch}{1.35}
\begin{tabular}{@{}llllc@{}}
\toprule
piece & range & defect $\delta=1-h'^2$ & $\kappa$ &
 lower bound for $h^2u^3k_{\widehat g}$\\
\midrule
cap, cylinder & $0\le r\le r_0$ & $\delta\equiv1$ \ ($h\equiv h_0$) & $1/6$ & $1/432$\\
ramp & $r_0\le r\le r_1$ & $hh''=\tfrac12\epsilon\eta\,\delta^2$ & $\delta/6$ & $\delta^2/500$\\
neck & $h_0\hn_*\le h\le e^{t_0}$ & $\dot\delta=-\epsilon\delta^2$ & $\delta/6$ & $\delta^2/500$\\
Stage I & $t_0\le t\le t_1$ & $\dot\delta=-\epsilon(1-\chi)\delta^2+\tfrac12\chi\delta$ & $\delta/6$ & $\delta^2/1000$\\
Stage II & $t_1\le t\le t_2$ & $\dot\delta=\tfrac12\delta$ & $\delta\psi/6$ & $\delta/12$\\
Stage III & $t_2\le t\le t_3$ & from $\delta_2e^{(t-t_2)/2}$ to $Ke^{2t}$ & $0$ & $\min\{\delta/2,\dot\delta/4\}$\\
\bottomrule
\end{tabular}
\caption{The pieces of the handle. Here $t=\log h$, a dot denotes
$d/dt$, and $\eta,\chi,\psi$ are cutoff functions; the precise
definitions are given in Subsection~\ref{subsec:profile}, and the lower bounds are
proved in Subsection~\ref{subsec:positivity}.}
\label{tab:pieces}
\end{table}

\subsection{The twisted metric and its curvature operator}
\label{subsec:curvature}

Let $r$ be a radial coordinate, and consider the metric
\begin{equation}\label{surg:ansatz}
 g=dr^2+f(r)^2d\theta^2+
 h(r)^2\bigl\{d\vartheta^2+
       \sin^2\vartheta\,(d\phi+A(r)\,d\theta)^2\bigr\},
\end{equation}
where $f$ and $h$ are positive smooth functions and $A$ is a smooth
function of $r>0$; the coordinate $\theta$ has period $2\pi$, and
$(\vartheta,\phi)$ are the spherical coordinates of
Theorem~\ref{surg:handle}. Thus $r$ is the distance in the radial
direction and $h$ is the radius of the sphere fibers. A prime denotes
the derivative with respect to $r$. We put
\begin{equation}\label{surg:notation}
\begin{aligned}
 \kappa &= \frac{h^2A'}{f}, &\quad z &= \cos\vartheta, &\quad
 b &= \frac{\kappa\sin\vartheta}{h}, &\quad
 c &= -\frac{\kappa\cos\vartheta}{h^2}, \\
 \cU &= \frac{1-h'^2}{h^2}, &\quad \cV &= -\frac{h''}{h}, &\quad
 \cW &= -\frac{f'h'}{fh}, &\quad \cX &= -\frac{f''}{f}.
\end{aligned}
\end{equation}
If $A$ is constant, then $\kappa=0$ and the metric \eqref{surg:ansatz} is
locally a warped product, and the functions $\cU,\cV,\cW,\cX$ are its
sectional curvatures. Here $\cU$ is the curvature of the planes tangent to
the spheres, $\cX$ that of the planes tangent to the disk, and $\cV$ and
$\cW$ are the curvatures of the planes spanned by a direction tangent to
the sphere and by $\partial_r$ or $\partial_\theta$, respectively. The
functions $b$ and $c$ measure the curvature of the connection. It is more
convenient to prescribe $\kappa$ than $A$.

On the region where $r>0$ and $0<\vartheta<\pi$ we use the oriented
orthonormal coframe
\[
 \omega^1=h\,d\vartheta,\qquad
 \omega^2=h\sin\vartheta\,(d\phi+A\,d\theta),\qquad
 \omega^3=dr,\qquad \omega^4=f\,d\theta,
\]
and we write $\omega^{ij}=\omega^i\wedge\omega^j$.

\begin{lemma}\label{lem:curvature-forms}
Let $\Omega_{ij}=\sum_{k<l}R_{ij,kl}\,\omega^{kl}$ be the curvature forms
of the metric \eqref{surg:ansatz} with respect to the coframe
$\omega^1,\ldots,\omega^4$, so that the coefficient of $\omega^{ij}$ in
$\Omega_{ij}$ is the sectional curvature of the plane spanned by $e_i$
and $e_j$. Then
\begin{equation}\label{surg:curvature-forms}
\begin{aligned}
\Omega_{12} &= \cU\omega^{12} + c\,\omega^{34} - \frac{h'b}{2h}\,\omega^{14}, &\qquad
\Omega_{13} &= \cV\omega^{13} + \frac{c}{2}\,\omega^{24}, \\
\Omega_{14} &= \cW\omega^{14} - \frac{c}{2}\,\omega^{23} - \frac{h'b}{2h}\,\omega^{12}, &\qquad
\Omega_{23} &= \Bigl(\cV + \frac{b^2}{4}\Bigr)\omega^{23} - \frac{c}{2}\,\omega^{14} + \mathcal{T}\omega^{34}, \\
\Omega_{24} &= \Bigl(\cW + \frac{b^2}{4}\Bigr)\omega^{24} + \frac{c}{2}\,\omega^{13}, &\qquad
\Omega_{34} &= \Bigl(\cX - \frac{3b^2}{4}\Bigr)\omega^{34} + c\,\omega^{12} + \mathcal{T}\omega^{23},
\end{aligned}
\end{equation}
where $\mathcal T=\dfrac{h'}{h}\,b+\dfrac{b'}{2}$ and
$b'=\partial b/\partial r$.
\end{lemma}

\begin{proof}
We abbreviate
\[
 \sigma=\frac{h'}h,\qquad \tau=\frac{f'}f,
 \qquad \mu=\frac{\cot\vartheta}{h},
\]
so that $c=-\mu b$, \ $\sigma'+\sigma^2=h''/h$ and
$\tau'+\tau^2=f''/f$. Differentiating the coframe, we find
\begin{align*}
 d\omega^1=\sigma\omega^{31},\qquad
 d\omega^2=\sigma\omega^{32}+\mu\,\omega^{12}+b\,\omega^{34},\qquad
 d\omega^3=0,\qquad
 d\omega^4=\tau\omega^{34},
\end{align*}
and the coefficients have the differentials
\[
 d\sigma=\sigma'\omega^3,\qquad d\tau=\tau'\omega^3,\qquad
 d\mu=-(h^{-2}+\mu^2)\,\omega^1-\sigma\mu\,\omega^3,
 \qquad db=\mu b\,\omega^1+b'\omega^3 .
\]
Let $\gamma_{ij}=-\gamma_{ji}$ be the connection forms, characterized by
$d\omega^i=-\sum_j\gamma_{ij}\wedge\omega^j$. One checks directly that
the structure equations are solved by
\begin{align*}
 \gamma_{12}&=-\mu\,\omega^2,&
 \gamma_{13}&=\sigma\omega^1,&
 \gamma_{14}&=0,\\
 \gamma_{23}&=\sigma\omega^2+\frac b2\,\omega^4,&
 \gamma_{24}&=-\frac b2\,\omega^3,&
 \gamma_{34}&=-\tau\omega^4-\frac b2\,\omega^2 .
\end{align*}
The curvature forms are
$\Omega_{ij}=d\gamma_{ij}+\sum_k\gamma_{ik}\wedge\gamma_{kj}$. For the
round unit sphere this convention gives $\Omega_{12}=\omega^{12}$, in
accordance with Section~\ref{sec:preliminaries}. We carry out two of the six
computations in detail. First,
\begin{align*}
 d\gamma_{12}
 =-d\mu\wedge\omega^2-\mu\,d\omega^2
   =h^{-2}\omega^{12}-\mu b\,\omega^{34},\qquad
 \gamma_{13}\wedge\gamma_{32}
 =-\sigma^2\omega^{12}-\frac{\sigma b}{2}\,\omega^{14},
\end{align*}
and $\gamma_{14}=0$, which gives the formula for $\Omega_{12}$ because
$\cU=h^{-2}-\sigma^2$. Second,
\begin{align*}
 d\gamma_{34}
 =\Bigl(\cX-\frac{b^2}{2}\Bigr)\omega^{34}
   +c\,\omega^{12}+\frac{b'+\sigma b}{2}\,\omega^{23},\qquad
 \gamma_{32}\wedge\gamma_{24}
 =-\frac{b^2}{4}\,\omega^{34}+\frac{\sigma b}{2}\,\omega^{23},
\end{align*}
and $\gamma_{31}\wedge\gamma_{14}=0$, which gives the formula for
$\Omega_{34}$. The remaining four forms are obtained in the same way from
\begin{align*}
 d\gamma_{13}&=\cV\omega^{13}, &
 \gamma_{12}\wedge\gamma_{23}&=\tfrac c2\,\omega^{24},\\
 d\gamma_{14}&=0, &
 \gamma_{12}\wedge\gamma_{24}+\gamma_{13}\wedge\gamma_{34}
 &=-\tfrac c2\,\omega^{23}+\cW\omega^{14}-\tfrac{\sigma b}2\,\omega^{12},\\
 d\gamma_{24}&=\tfrac c2\,\omega^{13}, &
 \gamma_{23}\wedge\gamma_{34}&=\bigl(\cW+\tfrac{b^2}4\bigr)\omega^{24},
\end{align*}
and
\begin{align*}
 d\gamma_{23}=\cV\omega^{23}+\sigma\mu\,\omega^{12}-\tfrac c2\,\omega^{14}
  +\bigl(\sigma b+\tfrac{b'}2+\tfrac{\tau b}2\bigr)\omega^{34},\\
 \gamma_{21}\wedge\gamma_{13}+\gamma_{24}\wedge\gamma_{43}
 =-\sigma\mu\,\omega^{12}-\tfrac{\tau b}2\,\omega^{34}+\tfrac{b^2}4\,\omega^{23}.
\end{align*}
As a check, the coefficients in \eqref{surg:curvature-forms} satisfy
$R_{ij,kl}=R_{kl,ij}$ and the Bianchi identity
$R_{12,34}+R_{13,42}+R_{14,23}=c-\frac c2-\frac c2=0$.
\end{proof}

We now pass to the Hodge decomposition. Recall that $A_+$ and $A_-$
denote the diagonal blocks of the curvature operator in
\eqref{eq:curvature-blocks}.

\begin{lemma}\label{lem:hodge-blocks}
With respect to the orthonormal bases
\[
 \frac{\omega^{12}\pm\omega^{34}}{\sqrt2},\qquad
 \frac{\omega^{13}\mp\omega^{24}}{\sqrt2},\qquad
 \frac{\omega^{14}\pm\omega^{23}}{\sqrt2}
\]
of $\Lambda^\pm$, the blocks $A_\pm$ of the metric \eqref{surg:ansatz}
are given by
\begin{equation}\label{surg:hodge}
 A_\pm=
 \begin{pmatrix}
  a_\pm&0&D\\
  0&\beta_\pm&0\\
  D&0&\beta_\pm
 \end{pmatrix},\qquad
 \begin{cases}
 a_\pm=\dfrac{\cU+\cX}{2}-\dfrac{3b^2}{8}\pm c,\\[5pt]
 \beta_\pm=\dfrac{\cV+\cW}{2}+\dfrac{b^2}{8}\mp\dfrac c2,\\[5pt]
 D=\dfrac{\kappa'\sin\vartheta}{4h}.
 \end{cases}
\end{equation}
Consequently, with the \emph{gaps} $G_\pm=a_\pm-\beta_\pm$, we have
\begin{equation}\label{surg:minimum}
 k_g=\frac{\cV+\cW}{2}+\frac{\kappa^2(1-z^2)}{8h^2}-\cE,
 \qquad
 \cE=\frac12\sum_{\pm}\Bigl(\sqrt{\tfrac14G_\pm^2+D^2}-\tfrac12G_\pm\Bigr)\ge0,
\end{equation}
and if $G_+$ and $G_-$ are positive, then
\begin{equation}\label{surg:error}
 \cE\le\frac{D^2}{\min\{G_+,G_-\}} .
\end{equation}
\end{lemma}

\begin{proof}
The entries of $A_\pm$ are read off from \eqref{surg:curvature-forms}.
For example, the entry coupling the first and the third basis vector is
\begin{align*}
 \frac12\bigl(R_{12,14}+R_{34,23}\pm R_{12,23}\pm R_{34,14}\bigr)
 =\frac12\Bigl(-\frac{h'b}{2h}+\mathcal T\Bigr)
 =\frac14\Bigl(b'+\frac{h'}hb\Bigr)
 =\frac{\kappa'\sin\vartheta}{4h},
\end{align*}
where we used $b=\kappa\sin\vartheta/h$ in the last step. In the same
way one finds the diagonal entries and sees that the second basis vector
is not coupled to the others. Thus $A_\pm$ is the direct sum of the
$1\times1$ matrix $(\beta_\pm)$ and of the symmetric $2\times2$ matrix
with diagonal entries $a_\pm,\beta_\pm$ and off-diagonal entry $D$. The
smaller eigenvalue of the latter is
\[
 \frac{a_\pm+\beta_\pm}2-\sqrt{\tfrac14G_\pm^2+D^2}
 =\beta_\pm-\cE_\pm,
 \qquad \cE_\pm=\sqrt{\tfrac14G_\pm^2+D^2}-\tfrac12G_\pm\ge0 .
\]
Hence $\lambda_{\min}(A_\pm)=\beta_\pm-\cE_\pm$, and \eqref{surg:minimum}
follows from Lemma~\ref{lem:hodge-minimum}, since
$\frac12(\beta_++\beta_-)=\frac12(\cV+\cW)+\frac18b^2$ and
$b^2=\kappa^2(1-z^2)/h^2$. If $G_\pm>0$, then
\[
 \cE_\pm=\frac{2D^2}{G_\pm+\sqrt{G_\pm^2+4D^2}}\le\frac{D^2}{G_\pm},
\]
which gives \eqref{surg:error}.
\end{proof}

We note two features of these formulas. First, the coupling $D$ involves
only the derivative of $\kappa$. If $\kappa$ is constant, the blocks
$A_\pm$ are diagonal, even where $h$ varies. Second, as a check on the
normalization, the scalar curvature is
$s_g=4\operatorname{tr}A_\pm=2\cU+4\cV+4\cW+2\cX-\frac12b^2$.

If $\kappa=0$, that is, if $A$ is constant, both blocks are diagonal with
entries $\frac12(\cU+\cX)$, $\frac12(\cV+\cW)$, $\frac12(\cV+\cW)$, so that
\begin{equation}\label{surg:untwisted}
 k_g=\min\Bigl\{\frac{\cU+\cX}{2},\frac{\cV+\cW}{2}\Bigr\}
 \qquad\text{if }\kappa=0 .
\end{equation}
Since $\frac12(\cV+\cW)=-(fh')'/(2fh)$, this gives
$k_g\le-(fh')'/(2fh)$, and where $f$ is constant it gives
$k_g\le-h''/(2h)$. In \eqref{surg:minimum}, by contrast, the term
$\frac12(\cV+\cW)$ is accompanied by the nonnegative connection term, and
the estimate holds for all two-planes, not only for those adapted to the
product structure.

\subsection{The conformal correction}\label{subsec:conformal}

The connection term in \eqref{surg:minimum} vanishes at the poles
$z=\pm1$, and since $\cE\ge0$ we have
$k_g\le\frac12(\cV+\cW)=-(fh')'/(2fh)$ there. As we saw at the beginning
of this section, the right-hand side is negative somewhere on the handle,
so that the twist alone does not suffice. We correct this by a conformal
change. For the metric \eqref{surg:ansatz} we put
\begin{equation}\label{surg:u}
 \cB=\frac{\kappa^2}{24},\qquad u=1+\cB z^2,
 \qquad \widehat g=u^2g ,
\end{equation}
and for a function $\varphi$ of $r$ we write
\[
 L_r\varphi=\varphi''+\Bigl(\frac{f'}f+\frac{2h'}h\Bigr)\varphi' .
\]

\begin{lemma}\label{lem:conformal-identity}
The metric $\widehat g$ in \eqref{surg:u} satisfies
\begin{equation}\label{surg:conformal-bound}
 u^3 k_{\widehat g}=
 -\frac u2\Bigl(\frac{h''}h+\frac{f'h'}{fh}\Bigr)
 +\frac{\kappa^2}{12h^2}
 +\frac{\kappa^4z^2(1-z^2)}{192h^2}
 -u\cE-\frac12(L_r\cB)\,z^2,
\end{equation}
where $\cE$ is the function in \eqref{surg:minimum}.
\end{lemma}

\begin{proof}
By Lemma~\ref{lem:conformal},
\begin{equation}\label{surg:conformal}
 u^3 k_{\widehat g}=u\,k_g-\frac12\Delta_g u .
\end{equation}
The volume density of $g$ in the coordinates
$(r,\theta,\vartheta,\phi)$ is $fh^2\sin\vartheta$, and for a function
of $r$ and $\vartheta$ alone the gradient has no component in the
directions of $\theta$ and $\phi$. Therefore
$\Delta_g\varphi=L_r\varphi$ for a function of $r$, and the Laplacian of
a function on the sphere is $h^{-2}$ times its Laplacian on the unit
sphere. Since $\Delta_{S^2}(z^2)=2-6z^2$, we get
\[
 \Delta_gu=(L_r\cB)\,z^2+\frac{\kappa^2}{24h^2}\,(2-6z^2) .
\]
We insert this and \eqref{surg:minimum} into \eqref{surg:conformal}. The
terms that involve neither derivatives of $h$ and $f$ nor $\cE$ combine to
\begin{align*}
 \frac{u\kappa^2(1-z^2)}{8h^2}
       -\frac{\kappa^2(2-6z^2)}{48h^2}
 &=\frac{\kappa^2}{48h^2}
       \bigl(6-6z^2-2+6z^2\bigr)
       +\frac{\kappa^4z^2(1-z^2)}{192h^2}\\
 &=\frac{\kappa^2}{12h^2}
       +\frac{\kappa^4z^2(1-z^2)}{192h^2},
\end{align*}
and \eqref{surg:conformal-bound} follows.
\end{proof}

The cancellation of the terms $\pm6z^2$ in this computation is the reason
for the constant $24$ in \eqref{surg:u}. It is the only choice for which
the leading term $\kappa^2/(12h^2)$ does not depend on $z$.

\subsection{The logarithmic radius}\label{subsec:log}

Outside the cap the function $f$ will be constant and $h$ will be
increasing. It is then convenient to use the logarithmic radius and the
defect,
\begin{equation}\label{surg:log-coordinates}
 t=\log h,\qquad \delta=1-h'^2,
\end{equation}
and to regard $\delta$ and $\kappa$ as functions of $t$. A dot denotes
the derivative with respect to $t$.

\begin{lemma}\label{lem:log}
Suppose that $f$ is constant and that $0<h'<1$ on some interval. Then
\begin{equation}\label{surg:log-derivatives}
 hh''=-\frac{\dot\delta}{2},\qquad
 \kappa'=\frac{\sqrt{1-\delta}}h\,\dot\kappa ,
\end{equation}
the gaps of Lemma~\ref{lem:hodge-blocks} are given by
\begin{equation}\label{surg:log-gap}
 2h^2G_\pm=
 \delta-\frac{\dot\delta}{2}
 -\kappa^2(1-z^2)\mp3\kappa z,
\end{equation}
and
\begin{equation}\label{surg:log-laplacian}
 h^2L_r\cB=(1-\delta)\ddot\cB+
          \Bigl(1-\delta-\frac{\dot\delta}{2}\Bigr)\dot\cB .
\end{equation}
Consequently, \eqref{surg:conformal-bound} takes the form
\begin{equation}\label{surg:log-master}
 h^2u^3k_{\widehat g}
 =\frac{u\,\dot\delta}{4}+\frac{\kappa^2}{12}
  +\frac{\kappa^4z^2(1-z^2)}{192}-u\,h^2\cE-\frac12\,(h^2L_r\cB)\,z^2 .
\end{equation}
If, moreover, $0\le\kappa\le\frac32$ and
$\delta-\frac12\dot\delta-3\kappa>0$, then
\begin{equation}\label{surg:log-error}
 2h^2G_\pm\ge\delta-\frac{\dot\delta}{2}-3\kappa
 \qquad\text{and}\qquad
 \cE\le\frac{(1-\delta)\,\dot\kappa^2}
        {8h^2\bigl(\delta-\frac12\dot\delta-3\kappa\bigr)} .
\end{equation}
\end{lemma}

\begin{proof}
Since $dt/dr=h'/h$, we have $\dot\delta=-2h'h''\cdot h/h'=-2hh''$ and
$\kappa'=\dot\kappa\,h'/h$, which is \eqref{surg:log-derivatives}. For
constant $f$ we have $\cW=\cX=0$, and by \eqref{surg:hodge}
\[
 2h^2G_\pm=h^2(\cU-\cV)-h^2b^2\pm3h^2c
 =\delta+hh''-\kappa^2(1-z^2)\mp3\kappa z,
\]
which gives \eqref{surg:log-gap}. For a function $\varphi$ of $r$ the
chain rule gives
\[
 \varphi'=\frac{h'}h\dot\varphi,\qquad
 \varphi''=\frac{h'^2}{h^2}\ddot\varphi
       +\Bigl(\frac{h''}h-\frac{h'^2}{h^2}\Bigr)\dot\varphi,
\]
and hence
$h^2L_r\varphi=h'^2\ddot\varphi+(hh''+h'^2)\dot\varphi$, which is
\eqref{surg:log-laplacian}. Now \eqref{surg:log-master} is
\eqref{surg:conformal-bound} multiplied by $h^2$.
For the last assertion we use the elementary fact that
\begin{equation}\label{surg:max}
 \max_{-1\le z\le1}
       \bigl\{\kappa^2(1-z^2)+3\kappa|z|\bigr\}=3\kappa
 \qquad\text{for }0\le\kappa\le\tfrac32 ;
\end{equation}
indeed, the derivative of $\kappa^2(1-x^2)+3\kappa x$ with respect to
$x=|z|$ is $\kappa(3-2\kappa x)\ge0$ for $0\le x\le1$. Together with
\eqref{surg:log-gap} this gives the lower bound for $2h^2G_\pm$. By
\eqref{surg:error}, \eqref{surg:hodge} and \eqref{surg:log-derivatives},
\begin{align*}
 \cE\le\frac{D^2}{\min\{G_+,G_-\}}
  \le\frac{(1-\delta)\,\dot\kappa^2(1-z^2)}{16h^4}\cdot
       \frac{2h^2}{\delta-\frac12\dot\delta-3\kappa}
   \le\frac{(1-\delta)\,\dot\kappa^2}
           {8h^2\bigl(\delta-\frac12\dot\delta-3\kappa\bigr)}.
   &\qedhere
\end{align*}
\end{proof}

The identity \eqref{surg:log-master} makes the strategy quantitative. With $\kappa=\delta/6$ the gap
condition reads $\frac12(\delta-\dot\delta)>0$, and the first two terms
on the right-hand side are $\frac14u\dot\delta+\frac1{432}\delta^2$. If
$\dot\delta=-\epsilon\delta^2$ with $\epsilon$ well below $\frac1{108}$,
their sum is a positive multiple of $\delta^2$, and the remaining terms
turn out to be of smaller order. We shall take $\epsilon=\frac1{1000}$.

\subsection{The radial profiles}\label{subsec:profile}

We now define the functions $f$, $h$ and $\kappa$. The function $A$ is
then determined by
\begin{equation}\label{surg:connection-integral}
 A'=\frac{f\kappa}{h^2},\qquad A(0)=0,
\end{equation}
and the metric of Theorem~\ref{surg:handle} will be
$\widehat g=u^2g$ with $g$ as in \eqref{surg:ansatz} and $u$ as in
\eqref{surg:u}. The outer pieces of the handle are fixed first, and the
radius $h_0$ of the sphere at the center is then chosen so that the neck
arrives at the first outer piece with the right defect.

\step{Constants and cutoff functions.}
Throughout, we fix
\[
 \epsilon=\frac1{1000},\qquad \kappa_0=\frac16 ,
\]
and a smooth nondecreasing function $\zeta$ from $\R$ to $[0,1]$ that vanishes
on a neighborhood of $(-\infty,0]$, is equal to one on a neighborhood
of $[1,\infty)$, and satisfies
\begin{equation}\label{surg:cutoff-bounds}
 0\le\zeta'\le2,\qquad |\zeta''|\le64 .
\end{equation}
Such a function exists. Indeed, a suitable rescaling of the function
$\exp(-1/(1-x^2))$ on $(-1,1)$ is a nonnegative smooth function with
integral one, with support in $[-\frac18,\frac18]$ and bounded by $16$.
The convolution of this function with $2\cdot\mathbf 1_{[1/4,3/4]}$ has a
primitive that vanishes at $0$, and this primitive has all the required
properties, since the absolute value of its second derivative is at most
$32$.

\step{The outer parameters.}
The last three pieces occupy fixed intervals of the logarithmic radius.
We put
\begin{equation}\label{surg:outer-parameters}
 t_3=\log R,\qquad t_2=t_3-1,\qquad
 t_1=t_3-11,\qquad t_0=t_3-31,
 \qquad \delta_0=\frac{KR^2}{4}\,e^{-17}.
\end{equation}
The number $\delta_0$ will be the defect at $t=t_0$. On the intervals
$[t_0,t_1]$ and $[t_1,t_2]$, of lengths $20$ and $10$, the defect will
satisfy $\dot\delta\le\frac12\delta$, so that it grows at most by the
factor $e^{15}$; the remaining factor $e^{-2}$ in $\delta_0$ ensures that
$\delta(t_2)\le\frac14Ke^{2t_2}$, which is what we need on the last
interval.

\step{The normalized ramp.}
Let $\hn$ be the solution on $[0,1]$ of the initial value problem
\begin{equation}\label{surg:normalized-ramp}
 \hn_{ss}=\frac{\epsilon\,\zeta(s)}{2\hn}\,(1-\hn_s^2)^2,
 \qquad \hn(0)=1,\qquad \hn_s(0)=0 .
\end{equation}
The pair $(\hn,p)$ with $p=\hn_s$ solves the first-order system
$\hn_s=p$, $p_s=\epsilon\zeta(1-p^2)^2/(2\hn)$. The set $\{p=1\}$ is
invariant under this system, and since $p(0)=0$, uniqueness shows that
$\hn_s<1$ as long as the solution exists. Hence
$0\le\hn_{ss}\le\frac12\epsilon$, and therefore
\[
 0\le\hn_s(s)\le\frac{\epsilon s}{2},\qquad
 1\le\hn(s)\le1+\frac{\epsilon s^2}{4},
\]
so that the solution exists on all of $[0,1]$. We put
\[
 \hn_*=\hn(1),\qquad \delta_*=1-\hn_s(1)^2 .
\]
Then $1\le\hn_*\le1+\frac14\epsilon$ and
$1-\frac14\epsilon^2\le\delta_*<1$; the last inequality holds because
$\zeta$ is positive on an interval, so that $\hn_s(1)>0$.

\step{The radius at the center.}
We define
\begin{equation}\label{surg:h0-choice}
 h_0=\frac{e^{t_0}}{\hn_*}\,
 \exp\Bigl[-\frac1\epsilon\Bigl(\frac1{\delta_0}-\frac1{\delta_*}\Bigr)\Bigr].
\end{equation}
Since $\delta_0<\delta_*$, we have $0<h_0\hn_*<e^{t_0}$.

\step{The cap and the cylinder, $0\le r\le r_0$.}
Put $\varrho(s)=1-\zeta(s-1)$ for $s\ge0$ and
$J=\int_0^\infty\varrho(s)\,ds$. The function $\varrho$ is nonincreasing,
equal to one near $[0,1]$, and it vanishes on $[2,\infty)$. We define
\begin{equation}\label{surg:f}
 f(r)=\frac FJ\int_0^{Jr/F}\varrho(s)\,ds,
 \qquad r_c=\frac{2F}{J},\qquad r_0=r_c+\lcyl,
\end{equation}
where $\lcyl\ge0$ is a parameter that will be fixed in
Subsection~\ref{subsec:end-of-proof}. Thus $f(r)=r$ near $r=0$, \ $f''\le0$, and
$f=F$ for $r\ge r_c$. On $[0,r_0]$ we put
\[
 h\equiv h_0,\qquad \kappa\equiv\kappa_0 .
\]
We call $[0,r_c]$ the \emph{cap} and $[r_c,r_0]$ the \emph{cylinder}.

\step{The ramp and the neck, $r\ge r_0$ and $h\le e^{t_0}$.}
From now on $f=F$. For $r\ge r_0$ we let $h$ be the solution of
\begin{equation}\label{surg:start-ode}
 h''=\frac{\epsilon\,\eta(r)}{2h}\,(1-h'^2)^2,
 \qquad h(r_0)=h_0,\quad h'(r_0)=0,
 \qquad\text{where }\eta(r)=\zeta\Bigl(\frac{r-r_0}{h_0}\Bigr),
\end{equation}
and we put
\[
 \kappa=\frac{\delta}{6}=\frac{1-h'^2}{6} .
\]
As for \eqref{surg:normalized-ramp}, we have $0\le h'<1$ and $h''\ge0$,
and the solution exists for all $r\ge r_0$. On the \emph{ramp}
$r_0\le r\le r_1$, where $r_1=r_0+h_0$, the substitution
$s=(r-r_0)/h_0$, $h=h_0\hn(s)$ transforms \eqref{surg:start-ode} into
\eqref{surg:normalized-ramp}. Hence
\[
 h(r_1)=h_0\hn_*,\qquad \delta(r_1)=\delta_*,\qquad h'(r_1)>0 .
\]
For $r\ge r_1$ we have $\eta\equiv1$ and $h'>0$, so that we can use the
logarithmic radius. By \eqref{surg:log-derivatives}, the equation
\eqref{surg:start-ode} becomes
\begin{equation}\label{surg:log-neck}
 \dot\delta=-\epsilon\delta^2,
 \qquad\text{that is,}\qquad
 \frac1{\delta}=\frac1{\delta_*}
       +\epsilon\log\frac{h}{h_0\hn_*} .
\end{equation}
The \emph{neck} is the region $h_0\hn_*\le h\le e^{t_0}$. The choice
\eqref{surg:h0-choice} of $h_0$ means precisely that
\[
 \delta=\delta_0\qquad\text{at } h=e^{t_0} .
\]

\step{Stage~I, $t_0\le t\le t_1$.}
In this stage the sign of $\dot\delta$ is reversed. Let
$\chi(t)=\zeta\bigl((t-t_0)/20\bigr)$, and let $\delta$ be the solution of
\begin{equation}\label{surg:stage-one-ode}
 \dot\delta=-\epsilon(1-\chi)\,\delta^2+\frac\chi2\,\delta,
 \qquad \delta(t_0)=\delta_0,
 \qquad\text{and put }\kappa=\frac\delta6 .
\end{equation}

\step{Stage~II, $t_1\le t\le t_2$.}
In this stage the connection is switched off. We let $\delta$ continue
to solve $\dot\delta=\frac12\delta$, and we put
\begin{equation}\label{surg:stage-two-def}
 \kappa=\frac{\delta\,\psi}{6},
 \qquad\text{where }\psi(t)=1-\zeta\Bigl(\frac{t-t_1}{10}\Bigr).
\end{equation}

\step{Stage~III, $t_2\le t\le t_3$.}
In the last stage the defect is steered to its value on the round collar.
With $\delta_2=\delta(t_2)$ we put $\kappa=0$ and
\begin{equation}\label{surg:stage-three}
 \delta(t)=\bigl(1-\zeta(t-t_2)\bigr)\,\delta_2e^{(t-t_2)/2}
                  +\zeta(t-t_2)\,Ke^{2t} .
\end{equation}

\medskip
In Stages~I--III the function $h$ is recovered from $\delta$ by solving
\begin{equation}\label{surg:recover-r}
 \frac{dr}{dt}=\frac{e^t}{\sqrt{1-\delta(t)}},
\end{equation}
which requires $\delta<1$; this is part of the following lemma.

\begin{lemma}\label{lem:profile}
The functions defined above have the following properties.
\begin{enumerate}[label=\textup{(\roman*)}]
\item On the ramp and the neck, $0<\delta\le1$, and
$\delta'=-\epsilon\eta h'\delta^2/h$.
\item In Stage~I, $\delta>0$ and $\dot\delta\le\frac12\delta$. In
Stages~I and~II,
\begin{equation}\label{surg:outer-smallness}
 0<\delta(t)\le\delta_0e^{15}
       =\frac{KR^2}{4}\,e^{-2}<\frac1{100},
 \qquad\text{and in particular }\ \delta_2\le\tfrac14Ke^{2t_2} .
\end{equation}
\item In Stage~III, $\dot\delta>0$ and
$0<\delta(t)\le Ke^{2t}\le KR^2\le\frac1{100}$.
\item The profiles fit together smoothly, and the handle has finite
length. More precisely, the coordinate $r$ runs through an interval
$[0,\ell]$, the functions $f$, $h$, $\kappa$ and $A$ are smooth on
$(0,\ell]$, and $h$ increases from $h_0$ to $R$ on $[r_0,\ell]$.
\item Near $r=\ell$ we have $\kappa=0$ and $\delta=Kh^2$, so that $u=1$,
$A$ is constant, and
\[
 \widehat g=g=F^2d\theta^2+\frac{dh^2}{1-Kh^2}
 +h^2\bigl\{d\vartheta^2+\sin^2\vartheta\,(d\phi+A(\ell)\,d\theta)^2\bigr\}.
\]
\end{enumerate}
\end{lemma}

\begin{proof}
(i) Differentiating $\delta=1-h'^2$ and using \eqref{surg:start-ode}
gives the formula for $\delta'$. We have already seen that $0\le h'<1$.

(ii) The function $y=1/\delta$ solves the linear equation
$\dot y=\epsilon(1-\chi)-\frac12\chi y$ with $y(t_0)>0$, and
$\epsilon(1-\chi)\ge0$; hence $y$ stays positive and finite on
$[t_0,t_1]$. The inequality $\dot\delta\le\frac12\delta$ is clear from
\eqref{surg:stage-one-ode}, and it also holds in Stage~II. Integrating it
over $[t_0,t_2]$, an interval of length $30$, gives
\eqref{surg:outer-smallness}; recall that $KR^2\le\frac1{100}$ and that
$t_2=t_3-1$.

(iii) Write $P(t)=\delta_2e^{(t-t_2)/2}$, $Q(t)=Ke^{2t}$ and
$\xi(t)=\zeta(t-t_2)$. By (ii), $Q/P\ge4e^{3(t-t_2)/2}>1$ on
$[t_2,t_3]$. Hence $\delta=(1-\xi)P+\xi Q$ lies between $P$ and $Q$, and
\begin{align*}
 \dot\delta
 =\tfrac12(1-\xi)P+2\xi Q+\dot\xi\,(Q-P)>0 .
\end{align*}

(iv) Each cutoff function is constant near the ends of the interval on
which it is used. Near $r=r_0$ we have $\eta=0$, so that $h$ and $\kappa$
are constant there. Near $t=t_0$ we have $\chi=0$, so that
\eqref{surg:stage-one-ode} agrees with \eqref{surg:log-neck}. Near
$t=t_1$ we have $\chi=1$ and $\psi=1$, so that
$\dot\delta=\frac12\delta$ and $\kappa=\frac16\delta$ on both sides, and
near $t=t_2$ we have $\psi=0$ and $\xi=0$, so that
$\dot\delta=\frac12\delta$ and $\kappa=0$ on both sides. Thus $\delta$
and $\kappa$ are smooth functions of $t$ on $[\log(h_0\hn_*),t_3]$, and
they are smooth functions of $r$ because $t$ is a smooth increasing
function of $r$ for $r\ge r_1$ by \eqref{surg:recover-r} and
\eqref{surg:start-ode}.
The neck has finite length, since $\delta\le\delta_*<1$ there, so that its
length is
\[
 \int_{h_0\hn_*}^{e^{t_0}}\frac{dh}{\sqrt{1-\delta}}
 \le\frac{e^{t_0}-h_0\hn_*}{\sqrt{1-\delta_*}}<\infty .
\]
Stages~I--III have finite length by \eqref{surg:recover-r}, since
$\delta\le\frac1{100}$ there. Finally, $A$ is smooth by
\eqref{surg:connection-integral}.

(v) Near $t=t_3$ we have $\xi=1$, hence $\delta=Ke^{2t}=Kh^2$, that is,
$h'^2=1-Kh^2$ and $dr^2=dh^2/(1-Kh^2)$. Since $\kappa=0$ there, $A$ is
constant and $u=1$.
\end{proof}

\subsection{Positivity}\label{subsec:positivity}

We now prove that $k_{\widehat g}>0$ on each piece of the handle. All
estimates are made on the region where $r>0$ and $0<\vartheta<\pi$, where
the coframe of Subsection~\ref{subsec:curvature} is defined; the center of the
disk and the poles of the sphere are treated at the end. Note that
$0\le\kappa\le\frac16$ everywhere, so that
\begin{equation}\label{surg:u-bounds}
 1\le u\le1+\frac{\kappa^2}{24}\le\frac{865}{864} .
\end{equation}

\begin{lemma}[Cap and cylinder]\label{lem:cap}
For $0<r\le r_0$,
\begin{equation}\label{surg:cap-positive}
 u^3 k_{\widehat g}
 =\frac{\kappa_0^2}{12h_0^2}
    +\frac{\kappa_0^4z^2(1-z^2)}{192h_0^2}
 \ge\frac1{432h_0^2} .
\end{equation}
\end{lemma}

\begin{proof}
Here $h$ and $\kappa$ are constant, so that $\cV=\cW=0$, $D=0$,
$L_r\cB=0$ and $\cU=h_0^{-2}$, while $\cX=-f''/f\ge0$. By
\eqref{surg:hodge} and \eqref{surg:max},
\[
 2h_0^2G_\pm=1+h_0^2\cX-\kappa_0^2(1-z^2)\mp3\kappa_0z
 \ge1-3\kappa_0=\frac12 .
\]
Since the gaps are positive and $D=0$, we have $\cE=0$, and
\eqref{surg:conformal-bound} reduces to \eqref{surg:cap-positive}.
\end{proof}

Note that the bound does not depend on the length $\lcyl$ of the cylinder.

\begin{lemma}[Ramp and neck]\label{lem:neck}
On the ramp and on the neck,
\begin{equation}\label{surg:neck-positive}
 u^3 k_{\widehat g}>\frac{\delta^2}{500\,h^2} .
\end{equation}
\end{lemma}

\begin{proof}
We work in the coordinate $r$, because $h'$ vanishes at the beginning of
the ramp. Recall that $f=F$, that $\kappa=\frac16\delta$, hence
$\cB=\delta^2/864$, and that by \eqref{surg:start-ode} and
Lemma~\ref{lem:profile}(i)
\begin{equation}\label{surg:delta-r}
 hh''=\frac{\epsilon\eta\,\delta^2}{2},\qquad
 h\delta'=-\epsilon\eta h'\delta^2 .
\end{equation}

\step{The radial term.} We claim that $L_r\cB\le0$. Differentiating
$\cB=\delta^2/864$ and using \eqref{surg:delta-r}, we find
\begin{align*}
 \cB'&=\frac{\delta\delta'}{432}
      =-\frac{\epsilon\eta h'\delta^3}{432h},\\
 \cB''&=-\frac\epsilon{432}
       \left[\frac{\eta'h'\delta^3}{h}
        +\frac{\eta h''\delta^3}{h}
        +\frac{3\eta h'\delta^2\delta'}{h}
        -\frac{\eta h'^2\delta^3}{h^2}\right],
\end{align*}
and therefore, since $h^2L_r\cB=h^2\cB''+2hh'\cB'$,
\begin{align*}
 h^2L_r\cB&=-\frac\epsilon{432}
       \left[\eta'hh'\delta^3+\eta hh''\delta^3
              +3\eta hh'\delta^2\delta'+\eta h'^2\delta^3\right].
\end{align*}
Substituting \eqref{surg:delta-r} gives
\begin{equation}\label{surg:start-radial}
 h^2L_r\cB=\frac1{864}\left[
 -2\epsilon\eta'hh'\delta^3
 +\epsilon^2\eta^2\delta^4(6h'^2-\delta)
 -2\epsilon\eta h'^2\delta^3\right].
\end{equation}
The first term is nonpositive because $\eta'\ge0$ and $h'\ge0$. With
$h'^2=1-\delta$, the other two terms add up to
\[
 \epsilon\eta\,\delta^3
 \bigl[\epsilon\eta\,\delta(6-7\delta)-2(1-\delta)\bigr].
\]
This is nonpositive. Indeed, for $\delta\ge\frac67$ both summands in the bracket
are nonpositive, and for $\delta<\frac67$ we use
$\delta(6-7\delta)\le6(1-\delta)$ and $\epsilon\eta\le\frac13$. Note that
no upper bound for $\eta'$ is needed.

\step{The gaps and the error term.} By \eqref{surg:hodge},
\eqref{surg:max} and \eqref{surg:delta-r},
\[
 2h^2G_\pm=\delta+hh''-\kappa^2(1-z^2)\mp3\kappa z
 \ge\delta+\frac{\epsilon\eta\,\delta^2}{2}-3\kappa
 \ge\frac\delta2 .
\]
Since $\kappa'=\frac16\delta'$, we obtain from \eqref{surg:error},
\eqref{surg:hodge} and \eqref{surg:delta-r}
\[
 \cE\le\frac{D^2}{\min\{G_+,G_-\}}
 \le\frac{\epsilon^2\eta^2h'^2\delta^4(1-z^2)}{576\,h^4}\cdot\frac{4h^2}{\delta}
 \le\frac{\epsilon^2\eta^2\delta^3}{144\,h^2} .
\]

\step{Conclusion.} We insert these estimates into
\eqref{surg:conformal-bound}, discard the nonnegative terms
$\kappa^4z^2(1-z^2)/(192h^2)$ and $-\frac12(L_r\cB)z^2$, and use
\eqref{surg:u-bounds}. This gives
\begin{align*}
 u^3 k_{\widehat g}
 &\ge-\frac{u\,\epsilon\eta\,\delta^2}{4h^2}+\frac{\delta^2}{432h^2}
   -\frac{u\,\epsilon^2\eta^2\delta^3}{144h^2}
 \ge\frac{\delta^2}{h^2}
 \left[\frac1{432}-\frac{865}{864}\Bigl(\frac\epsilon4
                   +\frac{\epsilon^2}{144}\Bigr)\right]
 >\frac{\delta^2}{500h^2} .
\end{align*}
The last inequality holds because the bracket is larger than
$0.00206$.
\end{proof}

\begin{lemma}[Stage~I]\label{lem:stage-one}
In Stage~I,
\begin{equation}\label{surg:stage-one-positive}
 u^3 k_{\widehat g}>\frac{\delta^2}{1000\,h^2} .
\end{equation}
\end{lemma}

\begin{proof}
Put $\lambda=\dot\delta/\delta$. By \eqref{surg:stage-one-ode},
$\lambda=-\epsilon(1-\chi)\delta+\frac12\chi$, and hence
\begin{equation}\label{surg:s-bounds}
 -\epsilon\delta\le\lambda\le\frac12,
 \qquad
 \dot\lambda
 =\dot\chi\,\bigl(\tfrac12+\epsilon\delta\bigr)
       -\epsilon(1-\chi)\lambda\delta\le\frac1{16} .
\end{equation}
For the last inequality we use that $0\le\dot\chi\le\frac1{10}$ by
\eqref{surg:cutoff-bounds} and that $\delta\le\frac1{100}$ by
Lemma~\ref{lem:profile}(ii); the term $-\epsilon(1-\chi)\lambda\delta$ is
positive only if $\lambda<0$, and then it is at most
$\epsilon^2\delta^2$.

We apply Lemma~\ref{lem:log}. Since $\kappa=\frac16\delta$,
\[
 \delta-\frac{\dot\delta}{2}-3\kappa
 =\frac{\delta(1-\lambda)}2\ge\frac\delta4 ,
\]
so that \eqref{surg:log-error} gives $h^2\cE\le\dot\kappa^2/(2\delta)$. As
$\dot\kappa=\frac16\lambda\delta$, we obtain
\[
 h^2\cE\le\frac{\lambda^2\delta}{72} .
\]
For $\cB=\delta^2/864$ we have $\dot\cB=\lambda\delta^2/432$ and
$\ddot\cB=\delta^2(\dot\lambda+2\lambda^2)/432$, and
\eqref{surg:log-laplacian} becomes
\begin{equation}\label{surg:stage-one-radial}
 h^2L_r\cB=\frac{\delta^2}{432}
 \left[(1-\delta)(\dot\lambda+2\lambda^2+\lambda)
                    -\frac{\lambda^2\delta}{2}\right]
 \le\frac{\delta^2}{432}\max\{0,\dot\lambda+2\lambda^2+\lambda\} .
\end{equation}
Inserting these estimates into \eqref{surg:log-master} and discarding
the quartic term, we obtain
\begin{equation}\label{surg:stage-one-estimate}
 h^2u^3 k_{\widehat g}\ge
 \frac{u\lambda\delta}{4}+\frac{\delta^2}{432}
 -\frac{u\lambda^2\delta}{72}
 -\frac{\delta^2}{864}
       \max\{0,\dot\lambda+2\lambda^2+\lambda\} .
\end{equation}
We distinguish two cases. If $0\le\lambda\le\frac12$, then
\begin{align*}
 \frac{u\lambda\delta}{4}-\frac{u\lambda^2\delta}{72}
 =u\lambda\delta\Bigl(\frac14-\frac\lambda{72}\Bigr)\ge0
 \qquad\text{and}\qquad
 \dot\lambda+2\lambda^2+\lambda
 \le\frac1{16}+\frac12+\frac12=\frac{17}{16},
\end{align*}
so that \eqref{surg:stage-one-estimate} gives
\[
 h^2u^3k_{\widehat g}
 \ge\delta^2\Bigl(\frac1{432}-\frac{17}{13824}\Bigr)
 =\frac{5\delta^2}{4608}>\frac{\delta^2}{1000} .
\]
If $-\epsilon\delta\le\lambda<0$, then $2\lambda^2+\lambda\le0$, so that
the maximum in \eqref{surg:stage-one-estimate} is at most $\frac1{16}$.
Using $|\lambda|\le\epsilon\delta$, $\delta\le1$ and
\eqref{surg:u-bounds}, we obtain
\[
 h^2u^3k_{\widehat g}\ge
 \delta^2\left[
 \frac1{432}-\frac{865}{864}
       \Bigl(\frac\epsilon4+\frac{\epsilon^2}{72}\Bigr)
 -\frac1{13824}\right] >\frac{\delta^2}{1000} ,
\]
since the bracket is larger than $0.0019$.
\end{proof}

\begin{lemma}[Stage~II]\label{lem:stage-two}
In Stage~II,
\begin{equation}\label{surg:stage-two-positive}
 u^3 k_{\widehat g}>\frac{\delta}{12\,h^2} .
\end{equation}
\end{lemma}

\begin{proof}
By \eqref{surg:cutoff-bounds}, the function $\psi$ in
\eqref{surg:stage-two-def} satisfies $0\le\psi\le1$,
$|\dot\psi|\le\frac15$ and $|\ddot\psi|\le\frac{64}{100}$; we shall only
use that $|\dot\psi|,|\ddot\psi|\le1$. From $\dot\delta=\frac12\delta$
and $\kappa=\frac16\delta\psi$ we get
\[
 \dot\kappa=\frac\delta6\Bigl(\dot\psi+\frac\psi2\Bigr),
 \qquad
 \ddot\kappa=\frac\delta6
       \Bigl(\ddot\psi+\dot\psi+\frac\psi4\Bigr),
 \qquad\text{hence}\qquad
 |\dot\kappa|\le\frac\delta4,
 \quad |\ddot\kappa|\le\frac{3\delta}{8} .
\]
Since
$\delta-\frac12\dot\delta-3\kappa=\frac34\delta-\frac12\delta\psi\ge\frac14\delta$,
Lemma~\ref{lem:log} applies, and \eqref{surg:log-error} gives
\[
 h^2\cE\le\frac{\dot\kappa^2}{2\delta}\le\frac{\delta}{32} .
\]
For $\cB=\kappa^2/24$ the bounds above imply
\[
 |\dot\cB|=\frac{|\kappa\dot\kappa|}{12}
            \le\frac{\delta^2}{288},
 \qquad
 |\ddot\cB|\le
       \frac{\dot\kappa^2+|\kappa\ddot\kappa|}{12}
            \le\frac{\delta^2}{96} .
\]
The coefficients in \eqref{surg:log-laplacian} are $1-\delta$ and
$1-\frac54\delta$, and both lie in $[0,1]$ because $\delta\le\frac1{100}$
by \eqref{surg:outer-smallness}. Hence
$|h^2L_r\cB|\le\frac1{96}\delta^2+\frac1{288}\delta^2=\frac1{72}\delta^2$.
Since $\frac14u\dot\delta=\frac18u\delta$, the identity
\eqref{surg:log-master} yields
\begin{align*}
 h^2u^3 k_{\widehat g}
 \ge\frac{u\delta}{8}-\frac{u\delta}{32}-\frac{\delta^2}{144}
 =\frac{3u\delta}{32}-\frac{\delta^2}{144}
 >\frac{\delta}{12},
\end{align*}
where we used $u\ge1$ and $\delta\le\frac1{100}$.
\end{proof}

\begin{lemma}[Stage~III]\label{lem:stage-three}
In Stage~III we have $\widehat g=g$ and
\begin{equation}\label{surg:untwisted-minimum}
 k_g=\frac1{h^2}
       \min\Bigl\{\frac\delta2,\frac{\dot\delta}{4}\Bigr\}>0 .
\end{equation}
Near $t=t_3$ this gives $k_g=\frac12K$.
\end{lemma}

\begin{proof}
Here $\kappa=0$, so that $u=1$ and $g$ is a warped product with constant
$f$. By \eqref{surg:untwisted} and \eqref{surg:log-derivatives},
$k_g=\min\{\frac12\cU,\frac12\cV\}$ with $\cU=\delta/h^2$ and
$\cV=\dot\delta/(2h^2)$. Both are positive by
Lemma~\ref{lem:profile}(iii). Near $t_3$ we have $\delta=Kh^2$ and
$\dot\delta=2Kh^2$.
\end{proof}
\subsection{Proof of Theorem~\ref{surg:handle}}\label{subsec:end-of-proof}

\begin{proof}[{\bf Proof of Theorem~\ref{surg:handle}}]
It remains to fix the length $\lcyl$ of the cylinder, to verify that
$\widehat g$ is a smooth metric on $D^2\times S^2$, and to collect the
results.

\step{Holonomy and framing.}
By Lemma~\ref{lem:profile}(v), the coefficient of $d\theta$ in the
boundary collar is the total rotation $A(\ell)=\int_0^\ell
f\kappa h^{-2}\,dr$. Let $A_0$ be its value for $\lcyl=0$. For general
$\lcyl\ge0$ the cylinder contributes $F\lcyl\kappa_0/h_0^2=F\lcyl/(6h_0^2)$ to this
integral, while all the later pieces are merely translated in $r$, so
that their contribution does not change. Thus
$A(\ell)=A_0+F\lcyl/(6h_0^2)$. We now choose
\begin{equation}\label{surg:n-L}
 n=\nu+2\Bigl\lceil\frac{A_0-\alpha-\nu}{2}\Bigr\rceil,
 \qquad
 \lcyl=\frac{6h_0^2}{F}\,(\alpha+n-A_0) .
\end{equation}
Then $n\equiv\nu\pmod2$. Moreover, $0\le\alpha+n-A_0<2$, so that
$\lcyl\ge0$, and $A(\ell)=\alpha+n$. By Lemma~\ref{lem:profile}(v), the metric near
$r=\ell$ is \eqref{surg:outer-model}. We define $R_-$ as a radius such
that this holds for $R_-<h\le R$.

\step{Smoothness at the center of the disk and at the poles.}
Near $r=0$ we have
\[
 f=r,\qquad h=h_0,\qquad \kappa=\frac16,
 \qquad A=\frac{r^2}{12h_0^2} .
\]
In the Cartesian coordinates $x=r\cos\theta$, $y=r\sin\theta$ on the disk,
\[
 dr^2+f^2d\theta^2=dx^2+dy^2,
 \qquad A\,d\theta=\frac{x\,dy-y\,dx}{12h_0^2},
\]
and the latter is a smooth one-form. In terms of the coordinates $(z,w)$
on the sphere, the fiber part of \eqref{surg:ansatz} is
\[
 h^2\bigl(dz^2+|dw+iw\,A\,d\theta|^2\bigr)
 \qquad\text{on } z^2+|w|^2=1,
\]
which is smooth at the poles $w=0$, and also where a pole meets the
center of the disk. Hence $g$ is a smooth metric on $D^2\times S^2$. The
conformal factor $u=1+\kappa^2z^2/24$ is smooth and positive, because
$\kappa$ is constant near $r=0$ and smooth for $r>0$. Therefore
$\widehat g=u^2g$ is a smooth Riemannian metric on $D^2\times S^2$.

\step{Positivity.}
By Lemmas~\ref{lem:cap}--\ref{lem:stage-three}, $k_{\widehat g}>0$ on the
dense open set where $r>0$ and $0<\vartheta<\pi$. In fact, the lower
bounds proved there are bounded below by a positive constant, because
$h\le R$ and $u\le\frac{865}{864}$, and because the continuous functions
$\delta$ on $[r_0,\ell]$ and $\dot\delta$ on $[t_2,t_3]$ are positive on
compact intervals. Since $k_{\widehat g}$ is continuous on
$D^2\times S^2$, it is positive everywhere. This completes the proof of
Theorem~\ref{surg:handle}.
\end{proof}

\begin{remark}\label{rem:size}
We have made no attempt to optimize the construction, and the handle is
extremely thin. With $KR^2=\frac1{100}$, the defect at the end of the
neck is $\delta_0\approx10^{-10}$, and \eqref{surg:h0-choice} gives
$h_0\approx R\exp(-10^{13})$. Accordingly, the total rotation $A_0$ of
the connection and the integer $n$ in \eqref{surg:n-L} are very large.
This is harmless, since only the parity of $n$ has a topological meaning
and we do not claim any effective bound on the geometry of the handle.
\end{remark}

\section{Proofs of the main results}
\label{sec:realization}

\subsection{Surgery along closed geodesics}

We first combine the results of Sections~\ref{sec:preparation} and~\ref{sec:handle}.

\begin{proposition}\label{prop:geodesic-surgery}
Let $(M,g)$ be a closed oriented Riemannian four-manifold with $k_g>0$,
and let $\gamma_1,\ldots,\gamma_m$ be pairwise disjoint embedded closed
geodesics of $g$. Then, for any choice of normal framings, the manifold
obtained from $M$ by surgery along $\gamma_1,\ldots,\gamma_m$ admits a
metric of positive biorthogonal curvature. This metric can be chosen to
agree with $g$ outside any given neighborhood of the geodesics.
\end{proposition}

\begin{proof}
We choose pairwise disjoint neighborhoods of the geodesics, as small as
we wish, and apply Lemma~\ref{lem:exact-tubes} with $K=1$. This deforms
$g$ inside these neighborhoods, through metrics of positive
biorthogonal curvature, to a metric that is a round tube near each
$\gamma_j$. Let $2\pi F_j$ be the length of $\gamma_j$. By
Lemma~\ref{lem:collar}, we can choose $\alpha_j\in\R$ and coordinates
$(\theta,h,(z,w))$ near $\gamma_j$ in which the new metric is
\eqref{eq:round-tube-boundary-form} with $F=F_j$, $K=1$ and
$\alpha=\alpha_j$, and we choose $\alpha_j$, within its class modulo
integers, such that these coordinates induce the given framing of
$\gamma_j$. We then fix $R_j>0$ with $R_j^2\le\frac1{100}$ so small that
the level set $\{h=R_j\}$ lies in the domain of these coordinates.

Theorem~\ref{surg:handle}, applied with
$(F,K,R,\alpha,\nu)=(F_j,1,R_j,\alpha_j,0)$, provides a metric of
positive biorthogonal curvature on $D^2\times S^2$ whose boundary collar
is \eqref{surg:outer-model} with an even integer $n_j$. We remove the open
tubes $\{h<R_j\}$ and glue in these handles by means of the maps
$T_{n_j}$, as described in Remark~\ref{rem:exact-collar-gluing}. The
metric is smooth across the seams because it is given by the same
expression on both sides, and it has positive biorthogonal curvature
everywhere. Since $n_j$ is even, the resulting manifold is the result of
surgery along the $\gamma_j$ with the given framings.
\end{proof}

\subsection{Proof of Theorem~\ref{thm:realization}}

\begin{proof}[{\bf Proof of Theorem~\ref{thm:realization}}]
\step{Step 1. Necessity, and the choice of a presentation.}
The fundamental group of a closed smooth manifold is finitely presented,
because the manifold has a finite triangulation.

Conversely, let $\Gamma$ be a finitely presented group. By
Lemma~\ref{lem:primitive-presentation} and its proof, $\Gamma$ has a
presentation with $N$ generators and $m$ relators that are primitive in
the free group $F_N$ and pairwise nonconjugate, even after inversion; if
we start from a presentation with $p$ generators and $q$ relators, none
of which is trivial in $F_p$, then $N=p+q$ and $m=2q$.

\step{Step 2. Geodesics that represent the relators.}
If $N=0$, then $\Gamma$ is trivial and the round four-sphere has the
required properties. Suppose that $N\ge1$, and let
\[
 M_0=\#^{N}(S^1\times S^3)
\]
be the connected sum of $N$ copies of $S^1\times S^3$, so that
$\pi_1(M_0)\cong F_N$. The product of a circle and a round three-sphere
has constant biorthogonal curvature $\frac12$, as we saw in the proof of
Lemma~\ref{lem:exact-tubes}, and Bettiol's connected-sum
theorem~\cite[Proposition~3.1]{Bettiol17} gives a metric $g_0$ on $M_0$
with $k_{g_0}>0$. If $m=0$, we are done.

Otherwise we apply Lemma~\ref{lem:geodesic-representatives} to the free
homotopy classes of the relators. Its hypotheses are satisfied, since a
primitive element of a free group is nontrivial and not a proper power,
and the relators are pairwise nonconjugate even after inversion. We
obtain a metric of positive biorthogonal curvature on $M_0$ and pairwise
disjoint embedded closed geodesics $\gamma_1,\ldots,\gamma_m$ of this
metric that represent the relators.

\step{Step 3. Surgery.}
Let $M$ be the manifold obtained from $M_0$ by surgery along
$\gamma_1,\ldots,\gamma_m$, with arbitrary framings. It is a closed
connected oriented smooth four-manifold, it admits a metric of positive
biorthogonal curvature by Proposition~\ref{prop:geodesic-surgery}, and by
Lemma~\ref{lem:circle-surgery}
\[
 \pi_1(M)\cong F_N/\langle\!\langle
     [\gamma_1],\ldots,[\gamma_m]
                         \rangle\!\rangle\cong\Gamma .
\]
This proves the theorem.
\end{proof}

\subsection{Proof of Corollary~\ref{cor:topology}}

\begin{proof}[{\bf Proof of Corollary~\ref{cor:topology}}]
Let $p$ and $q$ be as in the corollary. If $p+q=0$, the round sphere
satisfies all the assertions. Otherwise, let $M$ be the manifold
constructed above; it is obtained from $M_0=\#^{p+q}(S^1\times S^3)$ by
surgery along $2q$ circles, and we are free to choose the framings. Since
$\chi(M_0)=2-2(p+q)$ and since replacing $S^1\times D^3$ by
$D^2\times S^2$ increases the Euler characteristic by two,
\begin{equation}\label{eq:euler-presentation}
 \chi(M)=2-2(p+q)+4q=2(1-p+q).
\end{equation}
To make $M$ stably parallelizable, we regard $M_0$ as the boundary of the
five-dimensional handlebody $V$ with one zero-handle and $p+q$
one-handles. Since $V$ is homotopy equivalent to a graph,
$H^2(V;\Z/2)=0$, and $V$ is spin; we fix a spin structure on $V$. The
surgeries on $M_0=\partial V$ are the effect on the boundary of attaching
five-dimensional two-handles $D^2\times D^3$ to $V$ along the circles
$\gamma_j$. The two framings of $\gamma_j$ differ by the generator of
$\pi_1(\mathrm{SO}(3))\cong\Z/2$, and the spin structure of $V$ extends
over the corresponding two-handle for exactly one of
them~\cite[pp.~199--202]{MilnorSpin}. We choose these framings. The
resulting five-manifold $W$ is compact and spin, and $\partial W=M$.
Moreover, $W$ is homotopy equivalent to a two-dimensional complex, because
it is built from handles of index at most two. An oriented vector bundle
of rank five over a two-dimensional complex is trivial if its second
Stiefel--Whitney class vanishes. Indeed, $B\mathrm{SO}(5)$ is simply
connected and $\pi_2(B\mathrm{SO}(5))\cong\Z/2$, so that the only
obstruction to a null-homotopy of the classifying map lies in the second
cohomology group with coefficients in $\Z/2$, and this obstruction is the
second Stiefel--Whitney class. Hence the tangent bundle of $W$ is trivial.
Its restriction to $M$ is the sum of the tangent bundle of $M$ and a
trivial line bundle, and so $M$ is stably parallelizable.
Finally, Poincar\'e duality gives
$\chi(M)=2-2b_1(M)+b_2(M)$, and $b_1(M)=b_1(\Gamma)$ because
$H_1(M;\Q)$ is the rationalized abelianization of $\Gamma$. Together
with \eqref{eq:euler-presentation} this yields
$b_2(M)=2\bigl(q-p+b_1(\Gamma)\bigr)$.
\end{proof}

For later use we record what the proof gives.

\begin{remark}\label{rem:bounding}
The manifold $M$ of Corollary~\ref{cor:topology} is the boundary of a
compact parallelizable five-manifold $W$ such that the inclusion
$M\subset W$ induces an isomorphism of fundamental groups. Indeed, $W$ is
the union of $V$ and the trace of the surgeries, glued along $M_0$; the
inclusion $M_0\subset V$ induces an isomorphism of fundamental groups, and
so the claim follows from Lemma~\ref{lem:circle-surgery} and van Kampen's
theorem. By construction, $W$ has a handle decomposition with one
zero-handle, $p+q$ one-handles and $2q$ two-handles.
\end{remark}

\begin{remark}\label{rem:carr}
The auxiliary generators $y_j$ are needed to produce the geodesics, but
they do not affect the manifold. To see this, we slide, for each $j$, the
two-handle of $W$ that is attached along a circle representing $r_jy_j$
over the two-handle that is attached along a circle representing $y_j$.
For a suitable choice of the slides, the new attaching circles represent
$r_1,\ldots,r_q$. By Lemma~\ref{lem:isotopy} below we may assume that
these $q$ circles lie in the part of $M_0$ that comes from the first $p$
summands $S^1\times S^3$, and that the attaching circle representing
$y_j$ is the circle $S^1\times\{x\}$ in the summand that corresponds to
$y_j$. For each $j$, the one-handle that corresponds to $y_j$ and the
two-handle attached along $S^1\times\{x\}$ then cancel. Hence $W$ is
diffeomorphic to a parallelizable handlebody with one zero-handle, $p$
one-handles, and $q$ two-handles attached along circles that represent
$r_1,\ldots,r_q$. A regular neighborhood of the presentation complex of
$\langle x_1,\ldots,x_p\mid r_1,\ldots,r_q\rangle$, embedded in $\R^5$, is
a parallelizable handlebody of the same form. Thus the manifolds of
Corollary~\ref{cor:topology} are of the kind used by Carr~\cite{Carr88},
who proved that a regular neighborhood of a subcomplex of codimension at
least three in a Riemannian manifold can be chosen in such a way that its
boundary has positive scalar curvature. For positive scalar curvature,
statements corresponding to Corollary~\ref{cor:topology} and
Remark~\ref{rem:bounding} are recorded in~\cite[Theorem~6]{MT}.
\end{remark}

\begin{remark}\label{rem:group-examples}
Corollary~\ref{cor:topology} applies to groups with torsion as well. For
example, the presentation $\langle x\mid x^a\rangle$ of the cyclic group
of order $a\ge2$ has $p=q=1$ and gives a spin manifold with $\chi=2$ and
$b_2=0$, that is, a rational homology four-sphere. For the standard
presentation of $\Z^d$ with $d$ generators and $\binom d2$ commutator
relators we obtain
\[
 \chi(M)=(d-1)(d-2),\qquad b_2(M)=d(d-1).
\]
In particular, the manifold that realizes $\Z^4$ has $\chi=6$ and
$b_2=12$; it is not the four-torus, which admits no metric of positive
scalar curvature~\cite{SY79}. Theorem~\ref{thm:realization} prescribes
the fundamental group, not the manifold. We shall see in
Remark~\ref{rem:Z4} that the manifold constructed here has the smallest
Euler characteristic among all closed oriented four-manifolds of positive
scalar curvature with fundamental group $\Z^d$.
\end{remark}

\subsection{Proof of Theorem~\ref{thm:surgery}}

The proof of Theorem~\ref{thm:realization} uses surgery along closed
geodesics only. In order to pass to arbitrary circles we use two facts
from differential topology. The first is that in dimension four circles
cannot be knotted or linked.

\begin{lemma}\label{lem:isotopy}
Let $M$ be a closed connected smooth manifold of dimension at least four.
Let $c_1,\ldots,c_m$ and $c_1',\ldots,c_m'$ be two families of pairwise
disjoint embedded circles in $M$ such that $c_j$ is freely homotopic to
$c_j'$ for every $j$, for suitable orientations. Then there is a
diffeomorphism of $M$, isotopic to the identity, that maps $c_j$ onto
$c_j'$ for every $j$.
\end{lemma}

\begin{proof}
The two families define embeddings of the disjoint union of $m$ circles
into $M$, and these are homotopic by assumption. Since
$\dim M\ge 2\cdot1+2$, homotopic embeddings are isotopic by Whitney's
general position argument~\cite{Whitney36}; see also~\cite{Haefliger61}.
The reason is that during a generic one-parameter family of maps from a
one-dimensional manifold to a manifold of dimension at least four, no
double points and no singular points occur. By the isotopy extension
theorem~\cite[Chapter~8, Theorem~1.3]{Hirsch}, the isotopy is induced by
an isotopy of $M$.
\end{proof}

The second fact is that surgery along the circle factor undoes a
connected sum with $S^1\times S^3$.

\begin{lemma}\label{lem:cancel}
Surgery on $S^1\times S^3$ along the circle $S^1\times\{x\}$, with either
normal framing, yields $S^4$.
\end{lemma}

\begin{proof}
Write $S^3$ as the union of two balls $D^3_+$ and $D^3_-$ with $x$ the
center of $D^3_+$. Then $S^1\times D^3_+$ is a tubular neighborhood of
the circle, and the result of the surgery is
$D^2\times S^2\cup_\varphi S^1\times D^3_-$, where the gluing
diffeomorphism $\varphi$ of $S^1\times S^2$ is the identity for one
framing and the map $T_1$ of \eqref{eq:normal-gauge-twist} for the other.
Since $T_1$ extends to a diffeomorphism of $S^1\times D^3_-$, given by
the same formula, both manifolds are diffeomorphic to
$D^2\times S^2\cup_{\mathrm{id}}S^1\times D^3=\partial(D^2\times D^3)=S^4$.
\end{proof}

\begin{proof}[{\bf Proof of Theorem~\ref{thm:surgery}}]
Surgery along finitely many pairwise disjoint circles can be performed one
circle at a time, and the result of each step is again a closed connected
oriented manifold. It therefore suffices to consider a single circle.
Let $M$ be as in the theorem, let $c\subset M$ be an embedded circle, and
let $M'$ be the result of surgery along $c$ with one of its two framings.
We write $\Gamma=\pi_1(M)$, and we let $r\in\Gamma$ be an element in the
conjugacy class determined by $c$, for some orientation of $c$.

\step{Step 1. Stabilization.}
Let
\[
 M_2=M\,\#\,(S^1\times S^3)\,\#\,(S^1\times S^3),
\]
where the connected sums are formed in two balls that are disjoint from
$c$. Then $\pi_1(M_2)\cong\Gamma*F_2$, where the free group $F_2$ is
generated by the classes $y_1,y_2$ of the circles
$a_i=S^1\times\{x_i\}$ in the two new summands. Let $c''\subset M_2$ be
an embedded circle, disjoint from $a_1$ and $a_2$, that represents
$ry_1y_2$; for instance, we can take a connected sum of $c$ with parallel
copies of $a_1$ and $a_2$ along two bands.

We claim that $M'$ is obtained from $M_2$ by surgery along $a_1$, $a_2$
and $c''$, for arbitrary framings of $a_1$ and $a_2$ and a suitable
framing of $c''$. Indeed, by Lemma~\ref{lem:cancel}, surgery along $a_1$
and $a_2$ turns $M_2$ into $M\#S^4\#S^4=M$, and this identification is
the identity outside the two balls. By Lemma~\ref{lem:circle-surgery}, it
induces on fundamental groups the homomorphism $\Gamma*F_2\to\Gamma$ that
is the identity on $\Gamma$ and kills $y_1$ and $y_2$. Therefore $c''$,
regarded as a circle in $M$, is freely homotopic to $c$. By
Lemma~\ref{lem:isotopy}, there is a diffeomorphism of $M$ that maps $c''$
onto $c$; it maps the two framings of $c''$ to the two framings of $c$.
This proves the claim.

\step{Step 2. Geodesics.}
The manifold $M_2$ admits a metric of positive biorthogonal curvature by
Bettiol's connected-sum theorem. Consider the homomorphism
$\Gamma*F_2\to\Z^2$ that is trivial on $\Gamma$ and sends $y_1,y_2$ to
the standard basis vectors $e_1,e_2$. It maps $y_1$, $y_2$ and
$ry_1y_2$ to $e_1$, $e_2$ and $e_1+e_2$. None of these vectors is
divisible by an integer $d\ge2$, they are pairwise different, and none of
them is the negative of another. Hence the three elements are
nontrivial, they are not proper powers, and they are pairwise
nonconjugate, even after inversion.
Lemma~\ref{lem:geodesic-representatives} therefore provides a metric of
positive biorthogonal curvature on $M_2$ together with pairwise disjoint
embedded closed geodesics $\gamma_1,\gamma_2,\gamma_0$ in the free
homotopy classes of $a_1,a_2,c''$.

\step{Step 3. Conclusion.}
By Lemma~\ref{lem:isotopy}, there is a diffeomorphism of $M_2$ that maps
$a_1,a_2,c''$ onto $\gamma_1,\gamma_2,\gamma_0$. Hence, by Step~1, $M'$
is diffeomorphic to a manifold obtained from $M_2$ by surgery along the
geodesics $\gamma_1,\gamma_2,\gamma_0$, for suitable framings, and this
manifold admits a metric of positive biorthogonal curvature by
Proposition~\ref{prop:geodesic-surgery}.
\end{proof}

\begin{remark}\label{rem:no-stabilization}
The stabilization in Step~1 is needed only if the class of $c$ is a proper
power. Suppose first that this class is nontrivial and not a proper power.
Then Lemma~\ref{lem:geodesic-representatives} provides a metric of
positive biorthogonal curvature on $M$ and an embedded closed geodesic
$\gamma$ of this metric that is freely homotopic to $c$. By
Lemma~\ref{lem:isotopy} there is a diffeomorphism $\varphi$ of $M$ with
$\varphi(\gamma)=c$. Transporting the metric with $\varphi$, we obtain a
metric of positive biorthogonal curvature for which $c$ itself is a closed
geodesic, and Proposition~\ref{prop:geodesic-surgery} applies directly.
Suppose now that $c$ is null-homotopic. Then $c$ is isotopic to a circle
in a coordinate ball, and the result of the surgery is the connected sum
of $M$ with $S^2\times S^2$ or with the nontrivial $S^2$-bundle over
$S^2$, depending on the framing. Both admit positive biorthogonal
curvature by the results of Bettiol~\cite{Bettiol14,Bettiol17}.
\end{remark}

\begin{remark}\label{rem:surgery-implies-realization}
Theorem~\ref{thm:realization} also follows from
Theorem~\ref{thm:surgery}, by performing the surgeries along the relators
one after another; the passage from $r$ to $ry_1y_2$, $y_1$, $y_2$ in
Step~1 is the geometric counterpart of
Lemma~\ref{lem:primitive-presentation}. We have given the direct proof
because it treats all relators simultaneously and exhibits the manifold
as an explicit boundary, as needed in Corollary~\ref{cor:topology} and
Remark~\ref{rem:bounding}.
\end{remark}

\subsection{Proof of Corollary~\ref{cor:handlebody}}

\begin{proof}[{\bf Proof of Corollary~\ref{cor:handlebody}}]
By the standard rearrangement and cancellation of handles we may assume
that the handles of $W$ are attached in the order of their indices and,
since $W$ is connected, that there is exactly one zero-handle. Let
$W_1\subset W$ be the union of the zero-handle and the one-handles. Since
$W$ is oriented, $\partial W_1$ is the connected sum of $N$ copies of
$S^1\times S^3$, where $N$ is the number of one-handles, to be read as
$S^4$ if $N=0$. This manifold admits a metric of positive biorthogonal
curvature, as we saw in the proof of Theorem~\ref{thm:realization}. The
two-handles are attached to $W_1$ along pairwise disjoint embedded circles
in $\partial W_1$, and attaching a two-handle $D^2\times D^3$ changes the
boundary by surgery along its attaching circle. Hence $\partial W$ is
obtained from $\partial W_1$ by surgery along finitely many pairwise
disjoint embedded circles, and the assertion follows from
Theorem~\ref{thm:surgery}.
\end{proof}

\subsection{Proof of Corollary~\ref{cor:examples}}

\begin{proof}[{\bf Proof of Corollary~\ref{cor:examples}}]
\step{Part \textup{(a)}.}
Let $\Sigma$ be a closed oriented surface. Smale~\cite{Smale59} proved
that every oriented $S^2$-bundle over $\Sigma$ is the sphere bundle of an
oriented vector bundle of rank three. Hence it is the boundary of the
associated disk bundle $W\to\Sigma$. Choose a handle
decomposition of $\Sigma$. Since the bundle is trivial over each handle,
the part of $W$ over a handle of index $j$ is a five-dimensional handle of
index $j$. Thus $W$ is a compact connected oriented five-manifold that has
a handle decomposition with handles of index at most two, and
$\partial W$ admits a metric of positive biorthogonal curvature by
Corollary~\ref{cor:handlebody}. For the trivial bundle over the sphere
this gives a new proof that $S^2\times S^2$ admits positive biorthogonal
curvature~\cite{Bettiol14}.

\step{Part \textup{(b)}.}
By Perelman's work~\cite{Perelman03}, a closed orientable three-manifold
of positive scalar curvature is a connected sum of spherical space forms
$S^3/G$ and of copies of $S^2\times S^1$. In other words, it is obtained
from a disjoint union of spherical space forms by finitely many surgeries
of index zero that preserve orientability. We argue by induction on the
number of these surgeries. For a spherical space form, the product of
the round metric on $S^3/G$ with a circle is locally isometric to
$S^3\times\R$ and has biorthogonal curvature $\frac12$.

For the inductive step, let $Y'$ be obtained from $Y$ by a surgery of
index zero, where $Y$ is a possibly disconnected closed oriented
three-manifold, and suppose that each component of $Y\times S^1$ admits
positive biorthogonal curvature. The trace of the surgery is
$Y\times[0,1]$ with a one-handle $D^1\times D^3$ attached, and its
product with $S^1$ is a cobordism from $Y\times S^1$ to $Y'\times S^1$
that is obtained by attaching $D^1\times D^3\times S^1$. Writing $S^1$ as
the union of two arcs $I_1$ and $I_2$, we see that
$D^1\times D^3\times I_1$ is a five-dimensional one-handle and that
$D^1\times D^3\times I_2=(D^1\times I_2)\times D^3$ is a two-handle.
Hence $Y'\times S^1$ is obtained from $Y\times S^1$ by a surgery of index
zero followed by surgery along an embedded circle. The first of these
replaces $Y\times S^1$ by a connected sum of two of its components, or
by the connected sum of one component with $S^1\times S^3$, and preserves
positive biorthogonal curvature by Bettiol's theorem; the second
preserves it by Theorem~\ref{thm:surgery}.
\step{Parts \textup{(c)} and \textup{(d)}.}
Let $X$ be a compact connected oriented four-manifold that has a handle
decomposition with handles of index at most two. Then $X\times[0,1]$ has a
handle decomposition with handles of the same indices, and its boundary,
after the corners have been rounded, is the double of $X$. Let $Y$ be a
closed connected oriented three-manifold, and let $Y_0$ be the complement
of an open ball in $Y$. A handle decomposition of $Y$ with exactly one
three-handle restricts to a handle decomposition of $Y_0$ with handles of
index at most two, and $Y_0\times D^2$ has a handle decomposition with
handles of the same indices. Its boundary is the spun manifold of $Y$. In
both cases Corollary~\ref{cor:handlebody} applies.
\end{proof}

\section{Euler characteristics, stabilization, and open questions}
\label{sec:discussion}

\subsection{The minimal Euler characteristic}\label{subsec:euler}

Let $\Gamma$ be a finitely presented group. Following Hausmann and
Weinberger~\cite{HW85}, we write $q(\Gamma)$ for the minimum of the Euler
characteristics of closed connected oriented smooth four-manifolds with
fundamental group $\Gamma$. Following Mantione and Torres~\cite{MT}, we
write $q_{>0}(\Gamma)$ for the minimum over those manifolds that admit a
metric of positive scalar curvature, and we let $q_\perp(\Gamma)$ be the
minimum over those that admit a metric of positive biorthogonal
curvature. By Theorem~\ref{thm:realization} such manifolds exist, and
their Euler characteristics are bounded below because
$\chi(M)=2-2b_1(\Gamma)+b_2(M)$; hence the three minima exist. Since
positive biorthogonal curvature implies positive scalar curvature,
\begin{equation}\label{eq:three-invariants}
 q(\Gamma)\le q_{>0}(\Gamma)\le q_\perp(\Gamma).
\end{equation}

For a field $\F$ we put
$b_j(\,\cdot\,;\F)=\dim_{\F}H_j(\,\cdot\,;\F)$, both for spaces and for
groups, so that $b_j(\,\cdot\,)=b_j(\,\cdot\,;\Q)$. Recall that
the \emph{deficiency} $\defic(\Gamma)$ is the maximum of $p-q$ over all
finite presentations of $\Gamma$ with $p$ generators and $q$ relators. It
exists because
\begin{equation}\label{eq:deficiency-bound}
 p-q\le b_1(\Gamma;\F)-b_2(\Gamma;\F)
\end{equation}
for every such presentation and every field $\F$. Indeed, the presentation
complex $C$ satisfies $1-p+q=\chi(C)=1-b_1(C;\F)+b_2(C;\F)$. Moreover,
$H_1(C;\F)$ is isomorphic to $H_1(\Gamma;\F)$, and $H_2(C;\F)$ maps onto
$H_2(\Gamma;\F)$, because a space of type $K(\Gamma,1)$ can be obtained
from $C$ by attaching cells of dimension at least three.

\begin{corollary}\label{cor:deficiency}
Every finitely presented group $\Gamma$ satisfies
\[
 2-2b_1(\Gamma)+b_2(\Gamma)
 \le q(\Gamma)\le q_{>0}(\Gamma)\le q_\perp(\Gamma)
 \le 2-2\defic(\Gamma).
\]
All these inequalities are equalities for the free group $F_d$ of rank
$d\ge0$, where the common value is $2-2d$, and for nontrivial finite
cyclic groups, where it is $2$.
\end{corollary}

\begin{proof}
The first inequality is due to Hausmann and Weinberger~\cite{HW85}. It
holds because $\chi(M)=2-2b_1(\Gamma)+b_2(M)$ for every closed connected
oriented four-manifold $M$ with fundamental group $\Gamma$, and because
$H_2(M;\Q)$ maps onto $H_2(\Gamma;\Q)$, a space of type $K(\Gamma,1)$
being obtained from $M$ by attaching cells of dimension at least three.
The next two inequalities are \eqref{eq:three-invariants}. For the last
one we choose a presentation that realizes the deficiency. None of its
relators is trivial in the free group on the generators, since omitting
such a relator would increase $p-q$. Corollary~\ref{cor:topology}
therefore gives a manifold of positive biorthogonal curvature with
$\chi=2-2\defic(\Gamma)$. For $F_d$, the presentation with $d$ generators
and no relators shows that $\defic(F_d)=d=b_1(F_d)$, and $b_2(F_d)=0$.
For a nontrivial finite cyclic group, the presentation with one generator
and one relator shows that $\defic\ge0$, and $b_1=b_2=0$. In both cases
the two outer bounds agree.
\end{proof}

The inequality $q(\Gamma)\le2-2\defic(\Gamma)$ is also due to Hausmann and
Weinberger~\cite{HW85}, and the inequality
$q_{>0}(\Gamma)\le2-2\defic(\Gamma)$ follows from the construction of
Carr; see~\cite[Theorem~6]{MT}. What Corollary~\ref{cor:deficiency} adds is
that the manifolds which realize this upper bound admit positive
biorthogonal curvature.

The lower bound can be improved for manifolds whose fundamental class
vanishes in the homology of the group. The following observation is
known; see~\cite[Lemma~3.11]{Hillman02} for groups with
$H_4(\Gamma;\F)=0$, and~\cite{KL} and~\cite[Proposition~3]{MT} for
$\F=\Q$.

\begin{lemma}\label{lem:isotropic}
Let $M$ be a closed connected oriented four-manifold with fundamental
group $\Gamma$, let $c$ be a classifying map $M\to B\Gamma$ of its
universal cover, and let $\F$ be a field. If $c_*[M]=0$ in
$H_4(\Gamma;\F)$, then $b_2(M;\F)\ge2\,b_2(\Gamma;\F)$, and hence
\[
 \chi(M)\ge2\bigl(1-b_1(\Gamma;\F)+b_2(\Gamma;\F)\bigr).
\]
\end{lemma}

\begin{proof}
We may assume that $B\Gamma$ is obtained from $M$ by attaching cells of
dimension at least three and that $c$ is the inclusion. Then
the homomorphism $c^*$ from $H^2(\Gamma;\F)$ to $H^2(M;\F)$ is injective, and its image has
dimension $b_2(\Gamma;\F)$. For $a,b\in H^2(\Gamma;\F)$ we have
\[
 \langle c^*a\cup c^*b,[M]\rangle=\langle a\cup b,c_*[M]\rangle=0 ,
\]
so that the cup product pairing of $H^2(M;\F)$ vanishes on the image of
$c^*$. This pairing is nondegenerate by Poincar\'e duality, and therefore
a subspace on which it vanishes has at most half the dimension of
$H^2(M;\F)$. This proves the first assertion. The second one follows,
because $\chi(M)=2-2b_1(M;\F)+b_2(M;\F)$ and
$b_1(M;\F)=b_1(\Gamma;\F)$.
\end{proof}

\begin{corollary}\label{cor:minimal}
Let $\Gamma$ be a finitely presented group, and suppose that there is a
field $\F$ such that $H_4(\Gamma;\F)=0$ and
$\defic(\Gamma)=b_1(\Gamma;\F)-b_2(\Gamma;\F)$. Then
\[
 q(\Gamma)=q_{>0}(\Gamma)=q_\perp(\Gamma)=2-2\defic(\Gamma).
\]
\end{corollary}

\begin{proof}
By Lemma~\ref{lem:isotropic}, every closed connected oriented
four-manifold with fundamental group $\Gamma$ has Euler characteristic at
least $2\bigl(1-b_1(\Gamma;\F)+b_2(\Gamma;\F)\bigr)=2-2\defic(\Gamma)$.
The assertion now follows from Corollary~\ref{cor:deficiency}.
\end{proof}

For $q_{>0}(\Gamma)$ and $\F=\Q$, an equality of this kind was proved by
Mantione and Torres~\cite[Theorem~B]{MT}. For groups that admit a finite
aspherical two-dimensional complex we recover a theorem of
Hillman~\cite{Hillman94}, who determined $q(\Gamma)$ for these groups. In
both cases we see that the minimum is attained by manifolds of positive
biorthogonal curvature. Recall that a connected CW complex is
\emph{aspherical} if its universal cover is contractible.

\begin{corollary}\label{cor:aspherical-complex}
Let $X$ be a finite connected aspherical two-dimensional CW complex with
fundamental group $\Gamma$. Then
\[
 q(\Gamma)=q_{>0}(\Gamma)=q_\perp(\Gamma)=2\chi(X).
\]
In particular, for a closed oriented surface $\Sigma$ of positive genus,
\[
 q_\perp\bigl(\pi_1(\Sigma)\bigr)=2\chi(\Sigma),
\]
and the minimum is attained by both oriented $S^2$-bundles over
$\Sigma$.
\end{corollary}

\begin{proof}
Since $X$ is a space of type $K(\Gamma,1)$ of dimension two, we have
$H_4(\Gamma;\Q)=0$ and $\chi(X)=1-b_1(\Gamma)+b_2(\Gamma)$. Collapsing a
maximal tree in the one-skeleton of $X$, we obtain a presentation of
$\Gamma$ with $p$ generators and $q$ relators such that
$\chi(X)=1-p+q$. Hence $\defic(\Gamma)\ge p-q=b_1(\Gamma)-b_2(\Gamma)$,
and equality holds by \eqref{eq:deficiency-bound}.
Corollary~\ref{cor:minimal} now gives the first assertion, since
$2-2\defic(\Gamma)=2\chi(X)$. A closed oriented surface $\Sigma$ of
positive genus is aspherical. An $S^2$-bundle over $\Sigma$ has Euler
characteristic $\chi(S^2)\chi(\Sigma)=2\chi(\Sigma)$ and admits a metric
of positive biorthogonal curvature by Corollary~\ref{cor:examples}(a).
\end{proof}

\begin{remark}\label{rem:three-manifold-groups}
Corollary~\ref{cor:minimal} also applies to groups of cohomological
dimension three. Let $Y$ be a closed connected orientable aspherical
three-manifold with fundamental group $\Gamma$. Then $H_4(\Gamma;\Q)=0$,
and $b_1(\Gamma)=b_2(\Gamma)$ by Poincar\'e duality. A Heegaard splitting
of $Y$ gives a presentation of $\Gamma$ with as many relators as
generators, so that $\defic(\Gamma)\ge0$, and
\eqref{eq:deficiency-bound} shows that $\defic(\Gamma)=0$. Hence
$q(\Gamma)=q_\perp(\Gamma)=2$. The minimum is attained by the spun
manifold of $Y$, which has fundamental group $\Gamma$ and Euler
characteristic $2$ and admits positive biorthogonal curvature by
Corollary~\ref{cor:examples}(d). For $\Gamma=\Z^3$ it is also attained by
the manifold of Remark~\ref{rem:group-examples}.
\end{remark}

\begin{remark}\label{rem:Z4}
The invariants $q$ and $q_{>0}$ do not agree in general. The four-torus
and Corollary~\ref{cor:deficiency} show that $q(\Z^4)=0$, whereas
\[
 q_{>0}(\Z^d)=q_\perp(\Z^d)=(d-1)(d-2)
 \qquad\text{for every }d\ge1 ,
\]
so that $q_{>0}(\Z^4)=q_\perp(\Z^4)=6$. Indeed, let $M$ be a closed
connected oriented four-manifold with fundamental group $\Z^d$ that
admits a metric of positive scalar curvature, and let $c$ be a
classifying map $M\to T^d$. The classes of the four-dimensional
coordinate tori form a basis of $H_4(T^d;\Q)$. If $c_*[M]$ were nonzero
in this group, then the composition of $c$ with the projection of $T^d$
onto one of these coordinate tori would therefore have nonzero degree.
This is impossible, because a closed oriented four-manifold that admits a
map of nonzero degree to the four-torus carries no metric of positive
scalar curvature~\cite{SY79}. Hence Lemma~\ref{lem:isotropic} applies
with $\F=\Q$ and gives
$\chi(M)\ge2\bigl(1-d+\binom d2\bigr)=(d-1)(d-2)$. On the other hand,
the manifold of Remark~\ref{rem:group-examples} that realizes $\Z^d$ has
this Euler characteristic and admits positive biorthogonal curvature. For
$q_{>0}$, results of this kind are due to Mantione and
Torres~\cite[Theorem~B]{MT}, and the numbers $q(\Z^d)$ were determined
by Kirk and Livingston~\cite{KL}.
\end{remark}

We do not know whether $q_{>0}$ and $q_\perp$ can differ.

\begin{question}\label{q:qperp}
Is $q_\perp(\Gamma)=q_{>0}(\Gamma)$ for every finitely presented group
$\Gamma$?
\end{question}

A negative answer would also answer Question~\ref{q:psc}, since a manifold
of positive scalar curvature that realizes $q_{>0}(\Gamma)$ admits no
metric of positive biorthogonal curvature if
$q_{>0}(\Gamma)<q_\perp(\Gamma)$.

\subsection{Positive biorthogonal curvature after stabilization}

We now prove Corollary~\ref{cor:stable}. The proof combines
Theorem~\ref{thm:realization} with Kreck's stable classification of
four-manifolds~\cite[Theorem~C]{Kreck99}. It states that two closed
four-manifolds with the same Euler characteristic and the same normal
$1$-type that admit bordant normal $1$-smoothings become diffeomorphic
after connected sum with $r$ copies of $S^2\times S^2$, for some $r$. We
recall
from~\cite[Proposition~2]{Kreck99} what this means for a closed
connected oriented four-manifold $N$ with fundamental group $\Gamma$. If
the universal cover of $N$ is not spin, the normal $1$-type of $N$ is
$B\Gamma\times B\mathrm{SO}$; a normal $1$-smoothing is an orientation
together with a map $N\to B\Gamma$ that induces an isomorphism of
fundamental groups, and two such smoothings are bordant if they
represent the same element of the oriented bordism group
$\Omega_4^{\mathrm{SO}}(B\Gamma)$. If $N$ is spin, the same holds with
$B\mathrm{Spin}$, spin structures, and
$\Omega_4^{\mathrm{Spin}}(B\Gamma)$.

\begin{proof}[{\bf Proof of Corollary~\ref{cor:stable}}]
Let $M_\Gamma$ be a manifold as in Corollary~\ref{cor:topology}, for some
finite presentation of $\Gamma$, and let $W$ be the parallelizable five-manifold of
Remark~\ref{rem:bounding}. Since the inclusion induces an isomorphism
$\pi_1(M_\Gamma)\cong\pi_1(W)$, a classifying map
$c_\Gamma$ from $M_\Gamma$ to $B\Gamma$ extends to
$W$. Hence $(M_\Gamma,c_\Gamma)$ is null-bordant in
$\Omega_4^{\mathrm{SO}}(B\Gamma)$, and also in
$\Omega_4^{\mathrm{Spin}}(B\Gamma)$ for the spin structure induced from
$W$.

\step{Case \textup{(b)}.}
By assumption, $(N,c)$ and $(M_\Gamma,c_\Gamma)$ are bordant normal
$1$-smoothings in $B\Gamma\times B\mathrm{Spin}$. Both manifolds have
signature zero, so that their Euler characteristics have the same
parity. Connected sum with $S^2\times S^2$ raises the Euler
characteristic by two and changes neither the normal $1$-type nor the
bordism class, because $S^2\times S^2$ bounds the spin manifold
$S^2\times D^3$. After adding copies of $S^2\times S^2$ to one of the two
manifolds we may therefore assume that their Euler characteristics agree,
and Kreck's theorem shows that
\[
 N\,\#\,r\,(S^2\times S^2)\;\cong\;M_\Gamma\,\#\,r'(S^2\times S^2)
\]
for some $r,r'\ge0$. The manifold on the right admits a metric of
positive biorthogonal curvature, because $M_\Gamma$ and $S^2\times S^2$
do and the condition is preserved under connected
sums~\cite[Proposition~3.1]{Bettiol17}. Adding further copies of
$S^2\times S^2$ on both sides, we obtain the assertion for all larger
values of $r$.

\step{Case \textup{(a)}.}
Let $P$ be a connected sum of $a$ copies of the complex projective plane
and $b$ copies of the complex projective plane with the opposite
orientation, where $a,b\ge1$ and $a-b=\sign(N)$, and put
$M_1=M_\Gamma\#P$. Since the Fubini--Study metric has positive sectional
curvature, $M_1$ admits a metric of positive biorthogonal curvature. Its
fundamental group is $\Gamma$, and its universal cover is not spin,
because it contains a punctured copy of $P$ and hence an embedded sphere
with odd self-intersection number. Let $c_1$ be the map from $M_1$ to $B\Gamma$
that is $c_\Gamma$ on the first summand and constant on $P$; it
induces an isomorphism of fundamental groups, and $(c_1)_*[M_1]=0$.

Since $\Omega_j^{\mathrm{SO}}=0$ for $j=1,2,3$ and
$\Omega_4^{\mathrm{SO}}\cong\Z$ via the signature, the Atiyah--Hirzebruch
spectral sequence gives an isomorphism
\[
 \Omega_4^{\mathrm{SO}}(B\Gamma)\;\cong\;\Z\oplus H_4(\Gamma;\Z),
 \qquad [M,f]\longmapsto\bigl(\sign(M),\,f_*[M]\bigr).
\]
By assumption, $(N,c)$ and $(M_1,c_1)$ have the same image
$(\sign(N),0)$. Hence they are bordant normal $1$-smoothings in
$B\Gamma\times B\mathrm{SO}$. Since $N$ and $M_1$ have the same
signature, their Euler characteristics have the same parity, and we
conclude as in case~(b).
\end{proof}

\begin{remark}\label{rem:stable-general}
By the same argument, the conclusion of Corollary~\ref{cor:stable}
holds whenever a normal $1$-smoothing of $N$ is bordant to a normal
$1$-smoothing of a manifold of positive biorthogonal curvature. This
includes the third case of Kreck's classification, in which $N$ is not
spin but its universal cover is spin. In this case there is a unique class
$w\in H^2(\Gamma;\Z/2)$ with $c^*w=w_2(N)$, and the normal $1$-type of $N$
is determined by the pair $(\Gamma,w)$. In the proof of
Corollary~\ref{cor:topology} we can choose the framings of the two-handles
in such a way that the second Stiefel--Whitney class of $W$ is the image
of $w$ under the homomorphism induced by a classifying map
$W\to B\Gamma$, because changing the framing of a two-handle changes
$w_2(W)$ by the class of the corresponding two-cell, and these classes
generate $H^2(W;\Z/2)$. The boundary of the resulting manifold $W$ admits
positive biorthogonal curvature by Corollary~\ref{cor:handlebody}, it has
the same normal $1$-type as $N$, and it bounds over this normal $1$-type.
Hence the connected sum of $N$ with sufficiently many copies of
$S^2\times S^2$ admits positive biorthogonal curvature if some normal
$1$-smoothing of $N$ is null-bordant. Since positive biorthogonal
curvature implies positive scalar curvature, the manifolds in
Corollary~\ref{cor:stable} also admit positive scalar curvature after
stabilization.
\end{remark}

\subsection{Open questions}\label{subsec:questions}

\step{Positive scalar versus positive biorthogonal curvature.}
Question~\ref{q:psc} asks for an obstruction to positive biorthogonal
curvature that is not an obstruction to positive scalar curvature. The
Bochner technique does not provide one. By
Lemma~\ref{lem:hodge-minimum}, the condition $k_g>0$ means that
$\lambda_{\min}(A_+)+\lambda_{\min}(A_-)>0$, so that at each point at
least one of the operators $A_+$ and $A_-$ is positive definite. Where
$A_+$ is positive definite, the curvature term $\frac13s_g-2W_+$ in the
Weitzenb\"ock formula for self-dual harmonic two-forms is positive
definite as well. However, which of the two operators is positive may
depend on the point, and no vanishing theorem results. Indeed, $S^2\times S^2$
admits positive biorthogonal curvature, and both its self-dual and its
anti-self-dual harmonic forms are nontrivial. Likewise, by
Theorem~\ref{thm:realization}, the minimal surface arguments that
restrict the fundamental groups of manifolds of positive isotropic
curvature have no analogue for biorthogonal curvature. By
Theorem~\ref{thm:surgery} and Corollary~\ref{cor:stable}, whatever
distinguishes the two conditions is invisible to surgery of codimension
at least three and, in the range of Corollary~\ref{cor:stable}, to
connected sums with $S^2\times S^2$.

\step{Relative surgery.}
Proposition~\ref{prop:geodesic-surgery} is a relative statement, in the
sense that the metric is changed only near the geodesics.
Theorem~\ref{thm:surgery} holds for arbitrary embedded circles, but only
up to diffeomorphism, because its proof moves the circle to a geodesic of
a different metric on a stabilization of the manifold. We do not know
whether, for a given metric of positive biorthogonal curvature and an
arbitrary embedded circle, the surgery can be performed by a modification
of the metric inside a prescribed neighborhood of the circle. The handle
of Theorem~\ref{surg:handle} allows for both framings and for arbitrary
normal holonomy; what is missing is a deformation of the given metric to
a round tube near a circle that is not a geodesic.

\step{Nonorientable manifolds.}
Orientability of the manifold enters our arguments only through the
triviality of the normal bundle of the surgery circle. The constructions
of Sections~\ref{sec:preparation} and~\ref{sec:handle} take place in a
tubular neighborhood of the circle, and
Proposition~\ref{prop:geodesic-surgery} remains valid, with the same
proof, for closed geodesics with trivial normal bundle in a closed
nonorientable four-manifold of positive biorthogonal curvature. We expect
that Theorems~\ref{thm:realization} and~\ref{thm:surgery} have
nonorientable versions as well, in which the orientation character of the
manifold is prescribed, as in the result of Mantione and Torres for
positive scalar curvature~\cite[Theorem~7]{MT}. We do not pursue this
here.

\Needspace{9\baselineskip}
\section*{Acknowledgements}
This work was partially supported by NSFC (Grant Nos.\ 12571025, 12131012
and 12301032), NSF of Jiangsu (Grant No.\ BK20230803), the Fundamental
Research Funds for the Central Universities (Grant No.\ 4007012402), and
the Zhishan Scholars Program of Southeast University (Grant No.\
2242024RCB0039).

\medskip
\noindent\textit{Use of generative AI.} ChatGPT was used for exploratory
computations, proof development, and editorial assistance.

\noindent
Zhiqi Chen \\
School of Mathematics and Statistics, Guangdong University of Technology,\\
Guangzhou 510520, P.R. China\\
\textit{E-mail:} \href{mailto:chenzhiqi@nankai.edu.cn}{chenzhiqi@nankai.edu.cn}

\medskip
\noindent
Hui Zhang \\
School of Mathematics, Southeast University, Nanjing 211189, P.R. China\\
\textit{E-mail:} \href{mailto:huizhang@mail.nankai.edu.cn}{huizhang@mail.nankai.edu.cn}

\begin{thebibliography}{99}
\raggedright

\bibitem{BarHanke}
C.~B\"ar and B.~Hanke, \emph{Local flexibility for open partial differential
relations}, Comm. Pure Appl. Math. \textbf{75} (2022), no.~6, 1377--1415.
\href{https://doi.org/10.1002/cpa.21982}{doi:10.1002/cpa.21982}.

\bibitem{B1987}
A.~L.~Besse,
\emph{Einstein Manifolds},
Ergebnisse der Mathematik und ihrer Grenzgebiete (3), vol.~10,
Springer-Verlag, Berlin, 1987.

\bibitem{Bettiol14}
R.~G. Bettiol, \emph{Positive biorthogonal curvature on $S^2\times S^2$},
Proc. Amer. Math. Soc. \textbf{142} (2014), no.~12, 4341--4353.
\href{https://doi.org/10.1090/S0002-9939-2014-12173-3}{doi:10.1090/S0002-9939-2014-12173-3}.

\bibitem{Bettiol17}
R.~G. Bettiol, \emph{Four-dimensional manifolds with positive biorthogonal
curvature}, Asian J. Math. \textbf{21} (2017), no.~2, 391--396.
\href{https://doi.org/10.4310/AJM.2017.v21.n2.a7}{doi:10.4310/AJM.2017.v21.n2.a7}.

\bibitem{Carr88}
R.~Carr, \emph{Construction of manifolds of positive scalar curvature},
Trans. Amer. Math. Soc. \textbf{307} (1988), no.~1, 63--74.
\href{https://doi.org/10.1090/S0002-9947-1988-0936805-7}{doi:10.1090/S0002-9947-1988-0936805-7}.

\bibitem{CTZ}
B.-L. Chen, S.-H. Tang, and X.-P. Zhu, \emph{Complete classification of
compact four-manifolds with positive isotropic curvature},
J. Differential Geom. \textbf{91} (2012), no.~1, 41--80.
\href{https://doi.org/10.4310/jdg/1343133700}{doi:10.4310/jdg/1343133700}.

\bibitem{CZprod}
Z.~Chen and H.~Zhang, \emph{Positive biorthogonal curvature on $S^2$-bundles over  closed  surfaces}, preprint, 2026.
\href{https://arxiv.org/abs/2609.08119}{arXiv:2609.08119 v2}.

\bibitem{Costa12}
E. Costa, \emph{A modified Yamabe invariant and a Hopf conjecture},
preprint, 2012.
\href{https://arxiv.org/abs/1207.7107}{arXiv:1207.7107}.

\bibitem{CR14}
E. Costa and E. Ribeiro Jr., \emph{Four-dimensional compact manifolds
with nonnegative biorthogonal curvature}, Michigan Math. J.
\textbf{63} (2014), no.~4, 747--761.
\href{https://doi.org/10.1307/mmj/1417799224}{doi:10.1307/mmj/1417799224}.

\bibitem{Fraser03}
A.~M.~Fraser, \emph{Fundamental groups of manifolds with positive
isotropic curvature}, Ann. of Math. (2) \textbf{158} (2003), no.~1,
345--354.
\href{https://doi.org/10.4007/annals.2003.158.345}{doi:10.4007/annals.2003.158.345}.

\bibitem{GL80}
M.~Gromov and H.~B.~Lawson, Jr.,
\emph{The classification of simply connected manifolds of positive scalar curvature},
Ann. of Math. (2) \textbf{111} (1980), no.~3, 423--434.
\href{https://doi.org/10.2307/1971103}{doi:10.2307/1971103}.

\bibitem{Haefliger61}
A.~Haefliger, \emph{Plongements diff\'erentiables de vari\'et\'es dans
vari\'et\'es}, Comment. Math. Helv. \textbf{36} (1961), 47--82.

\bibitem{Hatcher}
A.~Hatcher, \emph{Algebraic Topology},
Cambridge University Press, Cambridge, 2002.

\bibitem{HW85}
J.-C.~Hausmann and S.~Weinberger,
\emph{Caract\'eristiques d'Euler et groupes fondamentaux des vari\'et\'es
de dimension $4$},
Comment. Math. Helv. \textbf{60} (1985), no.~1, 139--144.
\href{https://doi.org/10.1007/BF02567405}{doi:10.1007/BF02567405}.

\bibitem{Hillman94}
J.~A.~Hillman, \emph{Minimal $4$-manifolds for groups of cohomological dimension $2$},
Proc. Edinburgh Math. Soc. (2) \textbf{37} (1994), no.~3, 455--462.
\href{https://doi.org/10.1017/S0013091500018903}{doi:10.1017/S0013091500018903}.

\bibitem{Hillman02}
J.~A.~Hillman, \emph{Four-manifolds, geometries and knots},
Geometry \& Topology Monographs, vol.~5, Geometry \& Topology Publications,
Coventry, 2002.

\bibitem{Hirsch}
M.~W.~Hirsch, \emph{Differential Topology},
Graduate Texts in Mathematics, vol.~33, Springer-Verlag, New York, 1976.

\bibitem{Hoelzel}
S. Hoelzel, \emph{Surgery stable curvature conditions},
Math. Ann. \textbf{365} (2016), no.~1--2, 13--47.
\href{https://doi.org/10.1007/s00208-015-1265-1}{doi:10.1007/s00208-015-1265-1}.

\bibitem{KL}
P.~Kirk and C.~Livingston, \emph{The Hausmann--Weinberger $4$-manifold
invariant of abelian groups}, Proc. Amer. Math. Soc. \textbf{133} (2005),
no.~5, 1537--1546.

\bibitem{Kreck99}
M.~Kreck, \emph{Surgery and duality},
Ann. of Math. (2) \textbf{149} (1999), no.~3, 707--754.
\href{https://doi.org/10.2307/121071}{doi:10.2307/121071}.



\bibitem{MT}
A.~Mantione and R.~Torres, \emph{Geography of $4$-manifolds with positive
scalar curvature}, Expo. Math. \textbf{39} (2021), no.~4, 566--582,
\href{https://doi.org/10.1016/j.exmath.2021.05.003}{doi:10.1016/j.exmath.2021.05.003}.

\bibitem{MM88}
M.~J.~Micallef and J.~D.~Moore, \emph{Minimal two-spheres and the topology
of manifolds with positive curvature on totally isotropic two-planes},
Ann. of Math. (2) \textbf{127} (1988), no.~1, 199--227.
\href{https://doi.org/10.2307/1971420}{doi:10.2307/1971420}.

\bibitem{MilnorSpin}
J.~Milnor, \emph{Spin structures on manifolds},
Enseign. Math. (2) \textbf{9} (1963), 198--203.
\href{https://doi.org/10.5169/seals-38784}{doi:10.5169/seals-38784}.

\bibitem{MilnorStasheff}
J.~W.~Milnor and J.~D.~Stasheff, \emph{Characteristic Classes},
Annals of Mathematics Studies, vol.~76,
Princeton University Press, Princeton, NJ, 1974.

\bibitem{MP}
F. Morgan and P. Pansu, \emph{A list of open problems in differential
geometry}, S\~ao Paulo J. Math. Sci. \textbf{15} (2021), 305--321.
\href{https://doi.org/10.1007/s40863-019-00141-8}{doi:10.1007/s40863-019-00141-8}.

\bibitem{Perelman03}
G.~Perelman, \emph{Ricci flow with surgery on three-manifolds},
preprint, 2003.
\href{https://arxiv.org/abs/math/0303109}{arXiv:math/0303109}.



\bibitem{Rademacher}
H.-B.~Rademacher,
\emph{Simple closed geodesics in dimensions $\geq 3$},
J. Fixed Point Theory Appl. \textbf{26} (2024), article~5.
\href{https://doi.org/10.1007/s11784-023-01092-6}{doi:10.1007/s11784-023-01092-6}.

\bibitem{SY79}
R.~Schoen and S.-T.~Yau,
\emph{On the structure of manifolds with positive scalar curvature},
Manuscripta Math. \textbf{28} (1979), no.~1--3, 159--183.
\href{https://doi.org/10.1007/BF01647970}{doi:10.1007/BF01647970}.

\bibitem{Seaman91}
W. Seaman, \emph{Orthogonally pinched curvature tensors and applications},
Math. Scand. \textbf{69} (1991), 5--14.
\href{https://doi.org/10.7146/math.scand.a-12365}{doi:10.7146/math.scand.a-12365}.

\bibitem{ST69}
I.~M.~Singer and J.~A.~Thorpe, \emph{The curvature of $4$-dimensional
Einstein spaces}, in Global Analysis (Papers in Honor of K.~Kodaira),
Univ. Tokyo Press, Tokyo, 1969, pp.~355--365.

\bibitem{Smale59}
S.~Smale, \emph{Diffeomorphisms of the $2$-sphere},
Proc. Amer. Math. Soc. \textbf{10} (1959), 621--626.
\bibitem{Whitney36}
H.~Whitney, \emph{Differentiable manifolds},
Ann. of Math. (2) \textbf{37} (1936), no.~3, 645--680.
\href{https://doi.org/10.2307/1968482}{doi:10.2307/1968482}.

\end{thebibliography}
\end{document}